\documentclass[a4paper,oneside,11pt]{article}
 
\usepackage[a4paper]{geometry}
\usepackage{aeguill}                 
\usepackage{graphicx}
\usepackage{caption}
\usepackage{amsmath,amsfonts,amssymb,amsthm}
\usepackage{mathabx}
\usepackage{pstricks,comment} 
\usepackage{pst-plot}
\usepackage{pstricks-add}
\usepackage{multido}
\usepackage{url}
\usepackage{epstopdf,enumerate}
\usepackage{srcltx}
\usepackage[utf8]{inputenc} 
\usepackage[T1]{fontenc} 
\usepackage[english]{babel} 
\usepackage{hyperref} 
\usepackage[all]{xy} 
\usepackage{titlesec}
\usepackage{mathabx}
\usepackage{tikz}
\usepackage{fancyhdr}
\usepackage[us,12hr]{datetime}

    \newtheoremstyle{TheoremNum}
        {\topsep}{\topsep}              
        {\itshape}                      
        {}                              
        {\bfseries}                     
        {.}                             
        { }                             
        {\thmname{#1}\thmnote{ \bfseries #3}}

{\theoremstyle {definition} \newtheorem {defi} {Définition}[section]}
{\theoremstyle {plain}  \newtheorem {theo} [defi] {Théorème}}
{\theoremstyle {plain}  \newtheorem {coro} [defi] {Corollaire}}
{\theoremstyle {plain} \newtheorem {prop} [defi] {Proposition}}
{\theoremstyle {plain} \newtheorem {lem}[defi] {Lemme}}
{\theoremstyle {plain} \newtheorem {rmq}[defi] {Remarque}}
{\theoremstyle {plain} }
{\theoremstyle{TheoremNum} }
{\theoremstyle{TheoremNum} }
{\theoremstyle{TheoremNum} }

\newcommand{\Bij}{\mathrm{Bij}}
\newcommand{\Aut}{\mathrm{Aut}}
\newcommand{\Out}{\mathrm{Out}}
\newcommand{\Inn}{\mathrm{Inn}}

\newcommand{\PO}{\mathbb{P}\mathcal{O}}
\newcommand{\PA}{\mathbb{P}\mathcal{A}}
\newcommand{\Stab}{\mathrm{Stab}}
\newcommand{\Fix}{\mathrm{Fix}}
\newcommand{\ZZ}{\mathbb{Z}}

\newcommand{\RR}{\mathbb{R}}
\newcommand{\FF}{\mathbb{F}}

\newcommand{\dem}{\noindent{\bf Démonstration. }}

\addto{\captionsenglish}{}
\addto{\captionsenglish}{}

\title{Automorphismes du groupe des automorphismes \\ 
d'un groupe de Coxeter universel}
\author{Yassine Guerch}
\date{6 février 2020}
\begin{document}
\maketitle
\renewcommand{\proofname}{Preuve}
\renewcommand*\labelenumi{(\theenumi)}

\begin{abstract}
Using the Guirardel-Levitt outer space of a free product, we prove that the outer automorphism group of the outer automorphism group of the universal Coxeter group of rank $n\geq 5$ is trivial, and that it is a cyclic group of order $2$ if $n=4$. In addition we prove that the outer automorphism group of the automorphism group of the universal Coxeter group of rank $n \geq 4$ is trivial.
  \footnote{{\bf Keywords:} Universal Coxeter group, Outer automorphism groups, Outer space, Graph of groups.~~ {\bf AMS codes: } 20F55, 20E36, 20F28, 20E08}
\end{abstract}

\section{Introduction}

Soit $n$ un entier plus grand que $2$. On note $F= \ZZ/2\ZZ$ le groupe cyclique d'ordre $2$ et $W_n= \bigast_n F$ le groupe de Coxeter universel de rang $n$, produit libre de $n$ copies de $F$. Si $G$ est un groupe, on note $\Out(G)=\Aut(G)/\mathrm{Int}(G)$ son groupe d'automorphismes extérieurs. Nous démontrons dans cet article les résultats suivants.

\begin{theo}\label{Théorème out coxeter}
Si $n \geq 5$, alors $\Out\left(\Out\left(W_n\right)\right)=\{1\}$. Si $n=4$, alors $\Out(\Out(W_n))$ est isomorphe à $\ZZ/2\ZZ$.
\end{theo}

\begin{theo}\label{Théorème aut coxeter}
Si $n \geq 4$, alors $\Out(\Aut(W_n))=\{1\}$.
\end{theo}

De tels résultats sont déjà connus dans le cas où $n=2$ (cf. \cite[Lemma 1.4.2, Lemma 1.4.3]{thomasautotower}) où tous les automorphismes de $\Out(W_2)$ sont intérieurs et où $\Out(\Aut(W_2))$ est un groupe cyclique d'ordre $2$. Dans le cas où $n=3$, les groupes $\Aut(W_3)$ et $\Out(W_3)$ sont isomorphes à $\Aut(\FF_2)$ et $\mathbb{P}\mathrm{GL}(2,\ZZ)$ respectivement, avec $\FF_2$ un groupe libre de rang $2$ (cf. \cite[Lemma 2.3]{varghese2019}). Nous obtenons donc une description de $\Out(\Out(W_n))$ pour tout entier $n$ plus grand que $2$. 

\bigskip

De telles questions de rigidité algébrique ont déjà été résolues dans des cas similaires. En effet, Mostow \cite{Mostow73} a démontré que le groupe des automorphismes extérieurs de réseaux irréductibles uniformes de groupes de Lie réels, connexes, semi-simples et non localement isomorphes à $\mathrm{SL}_2(\RR)$ est fini. De même, Ivanov \cite[Theorem 2]{Ivanov97} a démontré un résultat similaire dans le cas du groupe modulaire d'une surface compacte, connexe, orientable de genre $g \geq 2$. Enfin, Bridson et Vogtmann \cite{bridson2000automorphisms} ont démontré que tout automorphisme du groupe des automorphismes extérieurs d'un groupe libre de rang $N$ (avec $N \geq 3$) est une conjugaison. Ce dernier cas a motivé l'étude de la rigidité algébrique de $\Out(W_n)$ d'une part à cause de la propriété d'universalité pour les groupes engendrés par des éléments d'ordre $2$ de $W_n$, d'autre part car, si $n \geq 3$, le groupe $\Aut(W_n)$ s'injecte dans $\Aut(F_{n-1})$ (cf. \cite[Theorem A]{muhlherr1997}). 

\bigskip

Pour démontrer les théorèmes~\ref{Théorème out coxeter} et \ref{Théorème aut coxeter}, nous étudions l'action de $W_n$ sur un complexe simplicial de drapeaux introduit par Guirardel et Levitt. Plus précisément, nous cherchons à comprendre les stabilisateurs de certains sommets de ce complexe. En effet, les stabilisateurs de ces sommets formant une partie génératrice de $\Aut(W_n)$ et $\Out(W_n)$, comprendre l'image de ces stabilisateurs par des automorphismes de $\Aut(W_n)$ et $\Out(W_n)$ nous permettra de faciliter l'étude de ces derniers. L'étude de l'action de $W_n$ sur un complexe simplicial se justifie également par la démonstration des théorèmes similaires dans les cas des réseaux des groupes de Lie semi-simples, du groupe modulaire d'une surface de type fini et du groupe des automorphismes d'un groupe libre qui passait également par l'étude de l'action du groupe étudié sur un espace géométrique adapté. En particulier, dans le cas du groupe des automorphismes extérieurs d'un groupe libre de rang $N$, cet objet géométrique était \emph{l'outre-espace de Culler-Vogtmann} $\mathrm{CV}_N$, qui fut introduit par Culler et Vogtmann dans \cite{Vogtmann1986}.

Dans le cas de $W_n$, Guirardel et Levitt \cite{Guirardel} ont introduit un espace topologique analogue à l'outre-espace de Culler et Vogtmann, appelé \emph{l'outre-espace d'un produit libre}. Dans le cas d'un produit libre de copies de $F$, cet espace sera noté $\PO(W_n)$. Ce dernier est défini comme un ensemble de classes d'homothétie de graphes de groupes métriques marqués de groupe fondamental isomorphe à $W_n$. Muni de la topologie dite \emph{faible}, l'espace $\PO(W_n)$ se rétracte par déformation forte sur un complexe simplicial de drapeaux, appelé \emph{l'épine de $\PO(W_n)$}. Le groupe $\Out(W_n)$ agit naturellement sur $\PO(W_n)$ et sur son épine par précomposition du marquage. Le groupe $\Aut(W_n)$ agit quant à lui sur \emph{l'autre espace de $W_n$}, noté $\PA(W_n)$. Nous renvoyons à la partie $2$ pour des précisions.

\bigskip

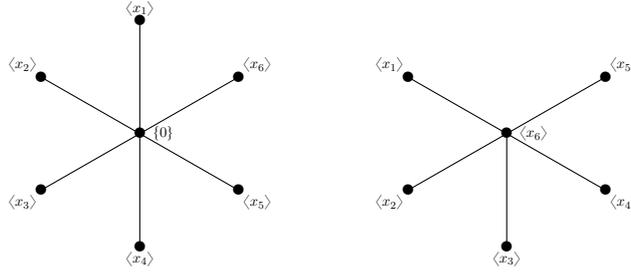
\begin{figure}[ht]
\centering
\captionsetup{justification=centering}
\begin{tikzpicture}[scale=1.5]
\draw (0:0) -- (90:1);
\draw (0:0) -- (150:1);
\draw (0:0) -- (210:1);
\draw (0:0) -- (270:1);
\draw (0:0) -- (330:1);
\draw (0:0) -- (30:1);
\draw (0:0) node {$\bullet$};
\draw (90:1) node {$\bullet$};
\draw (210:1) node {$\bullet$};
\draw (150:1) node {$\bullet$};
\draw (270:1) node {$\bullet$};
\draw (330:1) node {$\bullet$};
\draw (30:1) node {$\bullet$};
\draw (0:0) node[right, scale=0.5] {$\;\;\{0\}$};
\draw (90:1) node[above, scale=0.5] {$\left\langle x_1 \right\rangle$};
\draw (150:1) node[above left, scale=0.5] {$\left\langle x_2 \right\rangle$};
\draw (210:1) node[below left, scale=0.5] {$\left\langle x_3 \right\rangle$};
\draw (270:1) node[below, scale=0.5] {$\left\langle x_4 \right\rangle$};
\draw (330:1) node[below right, scale=0.5] {$\left\langle x_5 \right\rangle$};
\draw (30:1) node[above right, scale=0.5] {$\left\langle x_6 \right\rangle$};
\end{tikzpicture}
\hspace{1cm}
\begin{tikzpicture}[scale=1.5]
\draw (0:0) -- (270:1);
\draw (0:0) -- (150:1);
\draw (0:0) -- (210:1);
\draw (0:0) -- (330:1);
\draw (0:0) -- (30:1);
\draw (0:0) node {$\bullet$};
\draw (270:1) node {$\bullet$};
\draw (210:1) node {$\bullet$};
\draw (150:1) node {$\bullet$};
\draw (330:1) node {$\bullet$};
\draw (30:1) node {$\bullet$};
\draw (0:0) node[right, scale=0.5] {$\;\;\left\langle x_6 \right\rangle$};
\draw (270:1) node[below, scale=0.5] {$\left\langle x_3 \right\rangle$};
\draw (150:1) node[above left, scale=0.5] {$\left\langle x_1 \right\rangle$};
\draw (210:1) node[below left, scale=0.5] {$\left\langle x_2 \right\rangle$};

\draw (330:1) node[below right, scale=0.5] {$\left\langle x_4 \right\rangle$};
\draw (30:1) node[above right, scale=0.5] {$\left\langle x_5 \right\rangle$};
\end{tikzpicture}
\caption{Exemples de graphes de groupes dont les classes d'équivalence sont respectivement une $\{0\}$-étoile et une $F$-étoile (cas $n=6$). Les arêtes ont des groupes associés triviaux. L'ensemble $\{x_1,\ldots,x_6\}$ est une partie génératrice standard de $W_6$.}\label{figure 0étoile}
\end{figure}

La démonstration du théorème~\ref{Théorème out coxeter} est inspirée de celle de Bridson et Vogtmann dans le cas d'un groupe libre \cite{bridson2000automorphisms}, mais des complications structurelles apparaissent. Nous présentons la démonstration dans le cas de $\Out(W_n)$, le cas de $\Aut(W_n)$ étant similaire. Son plan, très simplifié, est le suivant. L'épine de l'outre-espace $\PO(W_n)$ contient, à la différence de celle de l'outre-espace de Culler-Vogtmann qui n'en contient qu'un, deux types de sommets distingués, à savoir les \emph{$\{0\}$-étoiles} et les \emph{$F$-étoiles}, voir la partie $2$ et la figure~\ref{figure 0étoile}. 

Nous étudions tout d'abord les stabilisateurs des $\{0\}$-étoiles et des $F$-étoiles sous l'action de $\Out(W_n)$. Nous montrons dans la partie $3$ que les sous-groupes de $\Out(W_n)$ isomorphes à $\mathfrak{S}_n$ sont les stabilisateurs de $\{0\}$-étoiles et les sous-groupes de $\Out(W_n)$ isomorphes au produit semi-direct $F^{n-2} \rtimes \mathfrak{S}_{n-1}$ sont les stabilisateurs de $F$-étoiles. Ces derniers représentent un cas nouveau en comparaison de la preuve de \cite{bridson2000automorphisms} dans le cas d'un groupe libre. De ce fait, tout automorphisme $\alpha$ de $\Out(W_n)$ préserve l'ensemble des stabilisateurs de $\{0\}$-étoiles et l'ensemble des stabilisateurs de $F$-étoiles. Fixons \mbox{$\alpha \in \Aut(\Out(W_n))$}. Le groupe $\Out(W_n)$ agissant transitivement sur l'ensemble des $\{0\}$-étoiles, nous pouvons supposer que $\alpha$ induit un automorphisme du stabilisateur d'une $\{0\}$-étoile $\mathcal{X}$. Les stabilisateurs de $\{0\}$-étoiles étant isomorphes à $\mathfrak{S}_n$, si $n \geq 5$ et $n \neq 6$, nous pouvons supposer que la restriction de $\alpha$ au stabilisateur de $\mathcal{X}$ est égale à l'identité. Nous montrons alors qu'un tel $\alpha$ préserve le stabilisateur d'une $F$-étoile $\mathcal{Y}$ adjacente à $\mathcal{X}$, et que la restriction de $\alpha$ au stabilisateur de $\mathcal{Y}$ est en fait l'identité. Le groupe $\Out(W_n)$ étant engendré par l'union des stabilisateurs d'une $\{0\}$-étoile et d'une $F$-étoile adjacente, ceci conclut la démonstration si $n \geq 5$. Le cas $n=4$, qui présente un automorphisme extérieur exceptionnel, est traité dans la partie $4$.

\bigskip

{\small{\bf Remerciements. } Je remercie chaleureusement mes directeurs de thèse, Camille Horbez et Frédéric Paulin, pour leurs précieux conseils et pour leur lecture attentive des différentes versions du présent article.}

\section{Préliminaires}
Nous rappelons tout d'abord la définition de l'outre-espace $\PO\left(W_n\right)$ introduit par Guirar\-del et Levitt dans \cite{Guirardel}. Un point de $\PO(W_n)$ est une classe d'homothétie de graphes de groupes métriques $X$ de groupe fondamental $W_n$ munis d'un isomorphisme de groupes appelé \emph{marquage} $\rho \colon W_n \to \pi_1(X)$ (pour un choix indifférent de point base) vérifiant~:
\begin{enumerate}
\item le graphe sous-jacent à $X$ est un arbre fini~;
\item tous les groupes d'arêtes sont triviaux~;
\item il y a exactement $n$ sommets de groupes associés isomorphes à $F$~;
\item tous les autres sommets ont un groupe associé trivial~;
\item toute feuille de l'arbre sous-jacent a un groupe associé non trivial~;
\item si $v$ est un sommet de groupe associé trivial, alors $\operatorname{deg}(v) \geq 3$.
\end{enumerate}
Deux graphes métriques marqués $(X,\rho)$ et $(X',\rho')$ sont dans la même classe d'homothétie s'il existe une homothétie $f \colon X \to X'$ (i.e.~un homéomorphisme multipliant toutes les longueurs des arêtes par un même scalaire strictement positif) telle que $f_{\ast} \circ \rho=\rho'$. On note $[X,\rho]$ la classe d'homothétie d'un tel graphe de groupes métrique marqué $(X,\rho)$. Si le marquage est sous-entendu, on notera $\mathcal{X}$ la classe d'homothétie. Le groupe $\Aut(W_n)$ agit par précomposition du marquage. Étant donné que pour tout $\alpha \in \Inn(W_n)$, et pour tout $\mathcal{X} \in \PO(W_n)$, nous avons $\alpha(\mathcal{X})=\mathcal{X}$, l'action de $\Aut(W_n)$ sur $\PO(W_n)$ induit une action de $\Out(W_n)$ sur $\PO(W_n)$.

La définition de \emph{l'autre espace de $W_n$}, noté $\PA(W_n)$, est identique à celle de $\PO(W_n)$ à ceci près que chaque graphe de groupes métrique considéré est muni d'un point base $v$. Le \emph{marquage} est alors un isomorphisme de groupes ${\rho \colon W_n \to \pi_1(X,v)}$. Les homothéties considérées préservent les points bases. Le groupe $\Aut(W_n)$ agit par précomposition du marquage.

\bigskip

L'ensemble $\PO(W_n)$ (resp. $\PA(W_n)$) est muni d'une topologie. Pour tout élément $[X,\rho] \in \PO(W_n)$, soit $(X,\rho)$ un représentant de cette classe d'équivalence tel que la somme des longueurs des arêtes du graphe sous-jacent soit égale à $1$. Le graphe de groupes $(X,\rho)$ définit alors un simplexe ouvert obtenu en faisant varier les longueurs des arêtes du graphe sous-jacent à $(X,\rho)$, de manière à ce que la somme des longueurs des arêtes soit toujours égale à $1$. Une classe d'équivalence $[X',\rho'] \in \PO(W_n)$ définit une face de codimension $1$ du simplexe associé à $(X,\rho)$ si l'on peut obtenir $(X',\rho')$ à partir de $(X,\rho)$ en écrasant une arête du graphe sous-jacent à $X$. La \emph{topologie faible} sur $\PO(W_n)$ est alors définie de la manière suivante~: un ensemble est ouvert si, et seulement si, son intersection avec chaque simplexe ouvert est ouverte. 

\bigskip
Nous rappelons à présent la définition d'un rétract par déformation forte $\Out(W_n)$-équivariant de $\PO(W_n)$, appelé \emph{l'épine de l'outre-espace}. L'épine de $\PO(W_n)$ est le complexe simplicial de drapeaux dont les sommets sont les simplexes ouverts associés à chaque classe d'équivalence $[X,\rho]$, et où deux sommets correspondant à des classes d'équivalence de graphes de groupes marqués $[X,\rho]$ et $[X',\rho']$ sont reliés par une arête si $[X,\rho]$ définit une face du simplexe associé à $[X',\rho']$ ou réciproquement. L'épine de $\PA(W_n)$ est définie de manière similaire. Il existe un plongement de l'épine de $\PO(W_n)$ dans $\PO(W_n)$ ayant pour image l'épine barycentrique de $\PO(W_n)$. Par la suite, nous identifierons l'épine de $\PO(W_n)$ avec son image par ce plongement. De même, il existe un plongement de l'épine de $\PA(W_n)$ dans $\PA(W_n)$ ayant pour image l'épine barycentrique de $\PA(W_n)$.
\bigskip

Si $X$ est un graphe de groupes, on note $\Aut_{gr}(X)$ le groupe des automorphismes du graphe sous-jacent à $X$. Si $X$ est un graphe de groupes pointé, la notation $\Aut_{gr}(X)$ désigne le groupe des automorphismes du graphe pointé sous-jacent à $X$.
Nous appellerons \emph{$\{0\}$-étoile} la classe d'équivalence dans $\PO(W_n)$ d'un graphe de groupes marqué dont le graphe sous-jacent est un arbre ayant $n$ feuilles et $n+1$ sommets et de longueur d'arêtes constante. Nous appellerons \emph{$F$-étoile} la classe d'équivalence dans $\PO(W_n)$ d'un graphe de groupes marqué dont le graphe sous-jacent est un arbre ayant $n-1$ feuilles et $n$ sommets et de longueur d'arêtes constante. Les sommets correspondants dans l'épine de $\PO(W_n)$ sont encore appelés \emph{$\{0\}$-étoiles} et \emph{$F$-étoiles}. Dans le cas de $\PA(W_n)$, les définitions des \emph{$\{0\}$-étoiles} et des \emph{$F$-étoiles} sont identiques à ceci près que l'on suppose également que le point base est le centre (l'unique sommet qui n'est pas une feuille) du graphe.

\bigskip

On fixe désormais une partie génératrice standard $\{x_1, \ldots,x_n\}$ de $W_n$.

Le groupe $\Aut\left(W_n\right)$ (et donc $\Out\left(W_n\right)$) est de type fini. Nous décrivons maintenant une partie génératrice finie. Pour tout $i \in \{1,\ldots,n-1\}$, on note $\tau_i \colon W_n \to W_n$ l'automorphisme envoyant $x_i$ sur $x_{i+1}$, $x_{i+1}$ sur $x_i$ et qui fixe tous les autres générateurs. Pour tous les $i,j \in \{1,\ldots,n\}$ tels que $i \neq j$, on note $\sigma_{i,j} \colon W_n \to W_n$ l'automorphisme qui envoie $x_i$ sur $x_jx_ix_j$ et qui fixe tous les autres générateurs. La proposition suivante est due à Mühlherr.

\begin{prop}\cite[Theorem B]{muhlherr1997}\label{système fini de générateurs}
Le groupe $\Aut(W_n)$ est engendré par $\tau_1,\ldots,\tau_{n-1}$ et par $\sigma_{1,2}$.
\end{prop}

Si $\alpha$ est un élément de $\Aut(W_n)$, sa classe d'automorphismes extérieurs sera notée $[\alpha]$. Soit $p \colon \Aut(W_n) \to \Out(W_n)$ la projection canonique.
On note $\widetilde{A}_n= \left\langle \tau_1,\ldots, \tau_{n-1} \right\rangle$ et $A_n=p(\widetilde{A}_n)$. Les groupes $\widetilde{A}_n$ et $A_n$ sont isomorphes au groupe symétrique $\mathfrak{S}_n$. On note $\widetilde{U}_n=\left\langle \tau_1,\ldots, \tau_{n-2}, \sigma_{1,n} \right\rangle$ et $U_n=p(\widetilde{U}_n)$. On voit que $\widetilde{U}_n$ est isomorphe au produit semi-direct $F^{n-1} \rtimes \mathfrak{S}_{n-1}$, alors que $U_n$ est isomorphe au produit semi-direct $F^{n-2} \rtimes \mathfrak{S}_{n-1}$, où $\mathfrak{S}_{n-1}$ agit dans les deux cas par permutation des facteurs, en considérant $F^{n-2}$ comme le quotient de $F^{n-1}$ par le sous-groupe $F$ diagonal. Soient $\widetilde{B}_n=\left\langle \tau_1,\ldots, \tau_{n-2} \right\rangle$ et $B_n=p(\widetilde{B}_n)$. Les groupes $\widetilde{B}_n$ et $B_n$ sont isomorphes à $\mathfrak{S}_{n-1}$.

\bigskip

Nous traitons à présent le cas où $n=3$. Soit $\epsilon \colon W_3 \to \ZZ/2\ZZ$ le morphisme envoyant, pour tout $i \in \{1,2,3\}$, l'élément $x_i$ sur $1$. Mühlherr (\cite[Theorem A]{muhlherr1997}) a montré que $\ker(\epsilon)$ est un sous-groupe caractéristique de $W_3$. De plus, $\ker(\epsilon)$ est un groupe libre à $2$ générateurs, librement engendré par $x_1x_2$ et $x_2x_3$. Ceci induit un morphisme $\rho \colon \Aut(W_3) \to \Aut(\FF_2)$, qui est en fait un isomorphisme (c.f. \cite[Lemma 2.3]{varghese2019}).

\begin{prop}
Le morphisme $\rho \colon \Aut(W_3) \to \Aut(\FF_2)$ induit un isomorphisme entre $\Out(W_3)$ et $\mathbb{P}\mathrm{GL}(2,\ZZ)$.
\end{prop}

\dem Soient $a$ et $b$ les générateurs de $\FF_2$. On remarque tout d'abord  que $\mathrm{Int}(\FF_2) \subseteq \rho(\mathrm{Int}(W_3))$. Donc le noyau du morphisme surjectif $\Aut(W_3) \to \Out(\FF_2)$ est inclus dans $\mathrm{Int}(W_3)$. Pour tout $i \in \{1,2,3\}$, soit $ad_{x_i} \in \Aut(W_3)$ la conjugaison globale par $x_i$. Un calcul immédiat montre que, pour tout $i \in \{1,2,3\}$, $\rho(ad_{x_i})$ est dans la classe d'automorphisme extérieur du morphisme $\iota \colon \FF_2 \to \FF_2$ envoyant $a$ sur $a^{-1}$ et $b$ sur $b^{-1}$. De ce fait, puisque le sous-groupe $\left\langle [\iota] \right\rangle$ est distingué dans $\Out(\FF_2)$, le morphisme $\rho$ induit un isomorphisme entre $\Out(W_3)$ et $\Out(\FF_2)/\left\langle [\iota] \right\rangle$. Comme $\iota$ est envoyé par le morphisme d'abélianisation sur $-\mathrm{Id} \in \mathrm{GL}(2,\ZZ)$, on voit que $\Out(W_3)$ est isomorphe à $\mathbb{P}\mathrm{GL}(2,\ZZ)$.
\hfill\qedsymbol

\bigskip

Nous allons démontrer les théorèmes~\ref{Théorème out coxeter} et \ref{Théorème aut coxeter} en étudiant les stabilisateurs des $\{0\}$-étoiles et des $F$-étoiles sous l'action de $\Out(W_n)$ et $\Aut(W_n)$. Pour cela, nous utiliserons les résultats suivants, dus respectivement à Hensel et Kielak et à Guirardel et Levitt, qui donnent des informations sur les points fixes de sous-groupes de $\Out(W_n)$.

\begin{prop}\cite[Corollary 6.1.]{hensel2018}\label{point fixe outre espace}
Soient $n \geq 1$ un entier et $H$ un sous-groupe fini de $\Out(W_n)$. Alors $H$ fixe un point de $\PO(W_n)$.
\end{prop}

\begin{coro}\label{points fixes autre espace aut}
Soient $n \geq 1$ un entier et $H$ un sous-groupe fini de $\Aut(W_n)$. Alors $H$ fixe un point de $\PA(W_n)$.
\end{coro}

\dem Soit $p \colon \Aut(W_n) \to \Out(W_n)$ la projection canonique. Alors $p(H)$ est un sous-groupe fini de $\Out(W_n)$, donc par la proposition~\ref{point fixe outre espace}, $p(H)$ fixe un point $\mathcal{X}$ de l'outre-espace. Soient $X$ un représentant de $\mathcal{X}$ et $\overline{X}$ son graphe sous-jacent. Comme tout automorphisme intérieur agit sur $X$, et que $p(H)$ agit également sur $X$, on en déduit que $H$ agit sur $X$. Étant donné que $H$ est fini et que $\overline{X}$ est un arbre, on voit que $H$ fixe un point $v$ de $\overline{X}$. Donc la classe d'homothétie du graphe de groupes métrique marqué pointé $(X,v)$ est fixée par $H$.
\hfill\qedsymbol

\begin{prop}\cite[Theorem 8.3.]{GuirardelDeformation}\label{Points fixes contractiles}
Soit $n \geq 2$ un entier. Si $H$ est un sous-groupe de type fini de $\Out(W_n)$ (resp. $\Aut(W_n)$) fixant un point de $\PO(W_n)$ (resp. $\PA(W_n)$), alors l'ensemble des points fixes de $H$ est contractile pour la topologie faible.
\hfill\qedsymbol
\end{prop}

On note $\mathrm{Fix}_{\PO(W_n)}(G)$ l'ensemble des points fixes d'un sous-groupe $G$ de $\Out(W_n)$ dans $\PO(W_n)$ (ou $\mathrm{Fix}(G)$ s'il n'y a pas d’ambiguïté). On note de plus $\mathrm{Fix}_{K_n}(G)$ l'ensemble des points fixes de $G$ contenus dans l'épine de $\PO(W_n)$.
Puisque l'épine de $\PO(W_n)$ est un rétract par déformation forte $\Out(W_n)$-équivariant de $\PO(W_n)$, nous déduisons de la proposition~\ref{Points fixes contractiles} le résultat suivant.

\begin{coro}\label{points fixes connexes}
Soit $n \geq 2$ un entier. Si $H$ est un sous-groupe de type fini de $\Out(W_n)$ fixant un point de l'épine de $\PO(W_n)$, alors l'ensemble $\Fix(H)$ des points fixes de $H$ dans l'épine de $\PO(W_n)$ est connexe pour la topologie faible.
\end{coro}

Soit $\mathcal{X}$ un point de l'épine de $\PO(W_n)$. On note $X$ un représentant de $\mathcal{X}$ et $T$ l'arbre de Bass-Serre associé à $X$. Nous définissons à présent un morphisme de groupes $$\Phi \colon \Stab_{\Out(W_n)}(\mathcal{X}) \to \Aut_{gr}(X).$$ Soient  $[\alpha] \in \Stab_{\Out(W_n)}(\mathcal{X})$, et $\alpha \in \Aut(W_n)$ un représentant de $[\alpha]$. Il existe un automorphisme $\widetilde{H}_{\alpha} \in \Aut (T)$ tel que pour tout $x \in T$, et pour tout $g \in W_n$ on ait \mbox{$\alpha(g) \widetilde{H}_{\alpha}(x)=\widetilde{H}_{\alpha}(gx)$}. L'automorphisme $\widetilde{H}_{\alpha}$ induit un automorphisme \mbox{$H_{\alpha} \in \Aut_{gr}(X)$}, et l'application $\alpha \mapsto H_{\alpha}$ passe au quotient pour donner un morphisme $$\Phi \colon \Stab_{\Out(W_n)}(\mathcal{X}) \to \Aut_{gr}(X).$$

\medskip

Nous pouvons à présent démontrer un résultat identique au corollaire~\ref{points fixes connexes} dans le cas de $\PA(W_n)$.

\begin{coro}\label{points fixes connexes Aut}
Soit $n \geq 2$ un entier. Si $H$ est un sous-groupe fini de $\Aut(W_n)$ fixant un point de l'épine de $\PA(W_n)$, alors l'ensemble $\Fix(H)$ des points fixes de $H$ dans l'épine de $\PA(W_n)$ est connexe pour la topologie faible.
\end{coro}

\dem
Soient $\mathcal{X}$ et $\mathcal{Y}$ deux points de l'épine de $\PA(W_n)$ fixés par $H$. Soit $p_1 \colon \PA(W_n) \to \PO(W_n)$ le morphisme canonique d'oubli du point base. On rappelle que \mbox{$p \colon \Aut(W_n) \to \Out(W_n)$} est la projection canonique. Alors $p(H)$ fixe $p_1(\mathcal{X})$ et $p_1(\mathcal{Y})$, donc par le corollaire~\ref{points fixes connexes} il existe dans $\mathrm{Fix}_{K_n}(p(H))$ un chemin continu $P$ de $p_1(\mathcal{X})$ vers $p_1(\mathcal{Y})$. Soient $\mathcal{X}_1,\ldots,\mathcal{X}_n$ les sommets de $K_n$ consécutifs dans $P$ (on suppose $p_1(\mathcal{X})=\mathcal{X}_1$ et $\mathcal{X}_n=p_1(\mathcal{Y})$) tels que, pour tout i $\in \{1,\ldots,n-1\}$, $X_i$ et $X_{i+1}$ sont reliés par une arête dans $K_n$. Soit $X_1$ un représentant de $\mathcal{X}_1$ et pour tout $i \in \{2,\ldots,n\}$, soit $X_i$ un représentant de $\mathcal{X}_i$ obtenu en écrasant ou en éclatant une forêt de $X_{i-1}$. Pour tout $i \in \{1,\ldots,n\}$, comme tout automorphisme intérieur agit trivialement sur $X_i$, et puisque $p(H)$ agit également sur $X_i$, on en déduit que $H$ agit sur $X_i$. De plus, étant donné que $H$ est fini et que le graphe sous-jacent $\overline{X}_i$ de $X_i$ est un arbre, on voit que $H$ fixe un point $v_i$ de $\overline{X}_i$. Pour tout $i$, soit $\widetilde{\mathcal{X}}_i$ la classe d'équivalence du graphe métrique marqué pointé $(X_i,v_i)$ (on suppose que $\widetilde{\mathcal{X}}_1=\mathcal{X}$ et $\widetilde{\mathcal{X}}_n=\mathcal{Y}$). Alors $\widetilde{\mathcal{X}}_i$ est fixé par $H$.

Nous construisons à présent pour tout $i \in \{1,\ldots,n-1\}$, un chemin continu inclus dans l'ensemble des points fixes de $H$ dans l'épine de $\PA(W_n)$ entre $\widetilde{\mathcal{X}}_i$ et $\widetilde{\mathcal{X}}_{i+1}$, ce qui conclura. La construction étant symétrique, nous pouvons supposer, quitte à changer les représentants $X_i$ et $X_{i+1}$, que $X_{i+1}$ est obtenu à partir de $X_i$ en écrasant une forêt $\mathcal{F}$. Soient $\Delta$ le simplexe ouvert dans $\PA(W_n)$ associé à $(X_i,v_i)$ et $e$ l'arête de l'épine barycentrique de $\PA(W_n)$ reliant $\widetilde{\mathcal{X}}_i$ et $\widetilde{\mathcal{X}}_{i+1}$. Pour toute arête $f$ de $\mathcal{F}$, soit $\ell_f$ la longueur de $f$. Pour tout $t \in [0,1]$, soient $X_i^t$ le graphe de groupes métrique obtenu à partir de $X_i$ en donnant à toute arête $f \in \mathcal{F}$ la longueur $(1-t)\ell_f$, et $pr_t \colon X_i \to X_i^t$ la projection canonique. On observe que $X_i^0=X_i$ et que $X_i^1=X_{i+1}$. 

Puisque $H$ stabilise $X_i$ et $X_{i+1}$, on voit que $H$ stabilise la forêt $F$. Donc, pour tout $t \in [0,1]$, le groupe $H$ stabilise $X_i^t$. Puisque $H$ fixe le sommet $v_i$ de $\overline{X}_i$, il fixe également, pour tout $t \in [0,1]$, le sommet $pr_t(x_i)$. Ceci induit un chemin continu de $\widetilde{\mathcal{X}}_i$ vers la classe d'équivalence dans $K_n$ de $(X_{i+1},pr_1(v_i))$. Si $pr_1(v_i) \neq v_{i+1}$, alors, puisque le graphe sous-jacent à $X_{i+1}$ est un arbre, $H$ fixe l'unique arc dans $\overline{X}_{i+1}$ reliant $pr_1(v_i)$ et $v_{i+1}$. Ceci induit alors un chemin continu contenu dans l'ensemble des points fixes de $H$ dans l'épine de $\PA(W_n)$ entre la classe d'équivalence dans $K_n$ de $(X_{i+1},pr_1(v_i))$ et $\widetilde{\mathcal{X}}_{i+1}$, ce qui conclut.
\hfill\qedsymbol

\bigskip

Soient $\mathcal{X}$ un point de l'épine de $\PO(W_n)$ et $X$ un représentant de $\mathcal{X}$. On note \mbox{$\Phi \colon \Stab_{\Out(W_n)}(\mathcal{X}) \to \Aut_{gr}(X)$} le morphisme naturel. Nous donnons maintenant une description de $\ker(\Phi)$. Soit $[X,\rho]$ un point de l'épine de $\PO(W_n)$. On note $(X,\rho)$ un représentant de $[X,\rho]$ et $\overline{X}$ le graphe sous-jacent à $X$. Soit $e$ une arête de $\overline{X}$ reliant le sommet $v=o(e)$ au sommet $w=t(e)$. Soit $z \in G_v$ un élément du groupe associé au sommet $v$, et $\overline{z}$ son antécédent par $\rho$. Nous définissons à présent le \emph{twist par $z$ autour de $e$}. Soit $G_u$ le groupe associé à un sommet $u$. Le twist par $z$ autour de $e$, noté $D_z$, est l'automorphisme de $W_n$, bien défini modulo conjugaison, qui est égal à l'identité sur $\rho^{-1}(G_u)$ si $u$ est dans la même composante connexe de $\overline{X}$ privé de l'intérieur de $e$ que $v$, et qui à $x \in \rho^{-1}(G_u)$ associe $\overline{z}x\overline{z}^{-1}$ si $u$ n'est pas dans la même composante connexe que $v$. Nous avons le résultat suivant, dû à Levitt.

\begin{prop}\cite[Proposition 2.2 and 3.1]{levitt2005}\label{prop Levitt}
Soit $n \geq 2$ un entier. Soient $\mathcal{X}$ un point de l'épine de l'outre-espace $\PO(W_n)$ et $X$ un représentant de $\mathcal{X}$. Soient $v_1,\ldots,v_n$ les sommets du graphe sous-jacent de $X$ de groupe associé isomorphe à $F$ et soit $n_i$ le degré de $v_i$ pour $i=1,\ldots,n$. Le noyau du morphisme $\Phi \colon \Stab_{\Out(W_n)}(\mathcal{X}) \to \Aut_{gr}(X)$ (noté $\Out_0(W_n)$ dans \cite{levitt2005}) est isomorphe à $\prod\limits_{i=1}^n F^{n_i-1}$, et il est engendré par les twists autour des arêtes dont l'origine appartient à $\{v_1,\ldots,v_n\}$ et n'est pas une feuille.
\end{prop}

\begin{rmq}\label{rmq twists Aut}
{\rm Dans le cas où $\mathcal{X} \in \PA(W_n)$, le noyau est engendré par les twists autour des arêtes $e$ dont l'origine $o(e)$ appartient à $\{v_1,\ldots,v_n\}$ et n'est pas une feuille, et telles que, si $o(e)$ est distinct du point base $v_{\ast}$, ces arêtes ne soient pas contenues dans l'unique chemin reliant $o(e)$ à $v_{\ast}$. En particulier, si le groupe associé à $v_{\ast}$ est trivial et si $n_i$ est le degré de $v_i$ pour $i=1,\ldots,n$, alors le noyau est isomorphe à $\prod\limits_{i=1}^{n} F^{n_i-1}$. Si le groupe associé à $v_{\ast}$ est non trivial, et si on suppose $v_{\ast}=v_n$, alors le noyau est isomorphe à $\left(\prod\limits_{i=1}^{n-1} F^{n_i-1}\right) \times F^{n_n}$.}
\end{rmq}

\section{Stabilisateurs des $\{0\}$-étoiles et des $F$-étoiles}

Nous étudions tout d'abord les stabilisateurs des $\{0\}$-étoiles.

\begin{lem}\label{Sn n feuilles}
Soit $n \geq 4$ un entier. Soient $G$ un sous-groupe fini de $\Out(W_n)$ isomorphe à $\mathfrak{S}_n$, et $\mathcal{X}$ un point de l'épine de $\PO(W_n)$ fixé par $G$. On note $X$ un représentant de $\mathcal{X}$ et $\overline{X}$ le graphe sous-jacent à $X$. Si le nombre de feuilles de $\overline{X}$ est $n$, alors $\mathcal{X}$ est une $\{0\}$-étoile.
\end{lem}

\dem Soit $v$ un sommet de $\overline{X}$ qui n'est pas une feuille et qui soit à distance maximale du centre\footnote{Rappelons que le centre d'un arbre métrique compact non vide est l'unique milieu d'un segment de longueur maximale.} de $\overline{X}$. 

\medskip

\noindent{\bf Affirmation. } Si $m=\operatorname{deg}(v)$, alors $v$ est adjacent à au moins $m-1$ feuilles de $\overline{X}$.

\medskip

\dem L'hypothèse de maximalité sur $v$ implique qu'il y a au plus un sommet adjacent à $v$ qui n'est pas une feuille, car sinon nous pourrions trouver un sommet $w$ adjacent à $v$ qui ne serait pas une feuille et qui serait à distance strictement plus grande du centre que $v$. \hfill\qedsymbol

\medskip

Maintenant, le groupe associé à $v$ est trivial car $\overline{X}$ possède exactement $n$ sommets de groupes associés non triviaux, et ces sommets sont tous des feuilles car $\overline{X}$ possède $n$ feuilles. De ce fait, $\operatorname{deg}(v) \geq 3$ et $v$ est adjacent à au moins deux feuilles, notées $v_1$ et $v_2$.

Soient $L$ l'ensemble des feuilles de $\overline{X}$, et $w$ une feuille de $\overline{X}$ distincte de $v_1$ et $v_2$. Puisque les seuls sommets de $\overline{X}$ dont les groupes associés sont non triviaux sont des feuilles, la proposition~\ref{prop Levitt} montre que le morphisme naturel $G \to \Aut_{gr}(X)$ est injectif. Ainsi, étant donné que le groupe $G$ est isomorphe à $\mathfrak{S}_n$, et que $\overline{X}$ possède $n$ feuilles, le morphisme naturel $\Aut_{gr}(X) \hookrightarrow \mathrm{Bij}(L)$ est un isomorphisme. Donc il existe un automorphisme de $\overline{X}$ envoyant $v_1$ sur $w$ et fixant $v_2$. De ce fait, $w$ est adjacent à $v$. Ainsi, $v$ est adjacent à toutes les feuilles de $\overline{X}$. Puisque le groupe $\Aut_{gr}(X)$ est isomorphe à $\Bij(L)$, toutes les arêtes de $\overline{X}$ ont même longueur. De ce fait, $\mathcal{X}$ est une $\{0\}$-étoile.
\hfill\qedsymbol

\begin{rmq}\label{Sn n feuilles aut}
{\rm Le résultat est identique dans le cas de $\PA(W_n)$. En effet, soit $G$ un sous-groupe fini de $\Aut(W_n)$ isomorphe à $\mathfrak{S}_n$, et $\mathcal{X}$ un point de l'épine de $\PA(W_n)$ fixé par $G$. On note $X$ un représentant de $\mathcal{X}$ et $\overline{X}$ le graphe sous-jacent à $X$. Supposons que $\overline{X}$ possède $n$ feuilles. Alors la remarque~\ref{rmq twists Aut} donne que le noyau du morphisme $G \to \Aut_{gr}(X)$ est un sous-groupe distingué de $G$ d'ordre au plus $2$. Comme $G$ est isomorphe à $\mathfrak{S}_n$, et que $n \geq 4$, le morphisme est injectif. La même démonstration que le lemme~\ref{Sn n feuilles} montre alors que $X$ possède $n$ feuilles et $n+1$ sommets. Il reste à montrer que le point base est le centre de $\overline{X}$.  Mais ceci provient du fait que le groupe $G$ est isomorphe à $\Aut_{gr}(X)$ qui lui-même est isomorphe à $\Bij(L)$. Ainsi, nécessairement, le point base est le centre de $\overline{X}$. Donc $\mathcal{X}$ est une $\{0\}$-étoile.}
\end{rmq}

\begin{prop}\label{Sn fixe 0 étoile}
Soient $n \geq 5$ un entier et $G$ un sous-groupe de $\Out(W_n)$ isomorphe à $\mathfrak{S}_n$. Alors $G$ est le stabilisateur dans l'épine de $\PO(W_n)$ d'une unique $\{0\}$-étoile.
\end{prop}

\dem Puisque $G$ est fini, d'après la proposition~\ref{point fixe outre espace}, il existe un point $\mathcal{X}$ de l'épine de l'outre-espace qui est fixé par $G$. Soit $X$ un représentant de $\mathcal{X}$. D'après la proposition~\ref{prop Levitt}, il existe un entier $k$ tel que le noyau de l'application naturelle \mbox{$G \to \Aut_{gr}(X)$} soit isomorphe à $F^k \cap G$.

Or $F^k \cap G$ est un $2$-sous-groupe distingué de $G \simeq \mathfrak{S}_n$. Donc, comme $n \geq 5$, un tel sous-groupe est trivial. De ce fait, $G$ s'injecte dans $\Aut_{gr}(X)$. Or tout automorphisme d'un arbre est entièrement déterminé par la permutation qu'il induit sur l'ensemble des feuilles. Ainsi, si $\overline{X}$ est le graphe sous-jacent à $X$ et si $L$ est l'ensemble des feuilles de $\overline{X}$,

\begin{equation*}
G \hookrightarrow \Aut_{gr}(X) \hookrightarrow \mathrm{Bij}(L).
\end{equation*}
Or les représentants des éléments de $\PO(W_n)$ possèdent au plus $n$ sommets de groupes non triviaux et toutes les feuilles possèdent des groupes associés non triviaux. Donc $|L| \leq n$. Donc, comme $G$ s'injecte dans $\mathrm{Bij}(L)$ et que $G$ est isomorphe à $\mathfrak{S}_n$, on voit que $G$ est isomorphe à $\Aut_{gr}(X)$ et que $\Aut_{gr}(X)$ est isomorphe à $\Bij(L)$. De ce fait, $\overline{X}$ possède $n$ feuilles. 
Par le lemme~\ref{Sn n feuilles}, $\mathcal{X}$ est une $\{0\}$-étoile.

Montrons maintenant l'unicité. Puisque l'ensemble des $\{0\}$-étoiles est discret dans l'épine de $\PO(W_n)$, par le corollaire~\ref{points fixes connexes}, on conclut que $G$ fixe une unique $\{0\}$-étoile dans l'épine de $\PO(W_n)$. 
\hfill\qedsymbol

\begin{rmq}\label{Sn fixe 0étoile aut}
{\rm Dans le cas de $\PA(W_n)$, le résultat de la proposition~\ref{Sn fixe 0 étoile} est vrai pour $n \geq 4$. En effet, dans le cas où $n \geq 5$, la démonstration est identique à celle de la proposition~\ref{Sn fixe 0 étoile} en utilisant cette fois la remarque~\ref{Sn n feuilles aut}.

\bigskip

Dans le cas où $n=4$, soit $\mathcal{X} \in \PA(W_4)$ un point fixé par un sous-groupe $G$ de $\Aut(W_n)$ isomorphe à $\mathfrak{S}_n$. On note $X$ un représentant de $\mathcal{X}$, $\overline{X}$ le graphe sous-jacent à $X$ et $v_{\ast}$ le point base de $\overline{X}$. Soit $H$ le noyau du morphisme $G \to \Aut_{gr}(X)$. Supposons par l'absurde que $H$ ne soit pas trivial. Alors, par la remarque~\ref{rmq twists Aut}, le groupe $H$ est un $2$-groupe. Comme le seul $2$-sous-groupe distingué de $\mathfrak{S}_4$ est le groupe de Klein, le groupe $H$ est isomorphe à $F^2$. Nous distinguons différents cas, selon le fait que le groupe associé à $v_{\ast}$ soit trivial ou non et selon le nombre de sommets qui ne sont pas des feuilles et qui ont un groupe associé non trivial. On remarque immédiatement que, puisque tout arbre possède au moins $2$ feuilles, le nombre de sommets qui ne sont pas des feuilles et de groupes associés non triviaux est au plus $2$.

\bigskip

Supposons que $\overline{X}$ contienne deux sommets qui ne soient pas des feuilles et dont les groupes associés sont isomorphes à $F$ et que le groupe associé à $v_{\ast}$ soit trivial. 

Soient $w_1$ et $w_2$ ces deux sommets. Alors $\operatorname{deg}(v_{\ast}) \geq 3$. Comme chaque composante connexe de $\overline{X}-\{v_{\ast}\}$ contient au moins une feuille, $\overline{X}$ contiendrait $5$ sommets de groupes associés non triviaux. Ceci contredit le fait qu'il y a exactement $4$ sommets dans le graphe de groupes associés non triviaux.

\bigskip
Supposons que $\overline{X}$ contienne deux sommets qui ne sont pas des feuilles et dont les groupes associés sont isomorphes à $F$ et que le groupe associé à $v_{\ast}$ ne soit pas trivial.

Alors la description du noyau du morphisme $G \to \Aut_{gr}(X)$ donné dans la remarque~\ref{rmq twists Aut} donne que le cardinal du noyau est au moins $8$, ce qui contredit le fait que $H$ est de cardinal $4$.

\bigskip

Supposons que $\overline{X}$ contienne un seul sommet, noté $w$, de groupe associé non trivial et qui ne soit pas une feuille et que le groupe associé à $v_{\ast}$ soit trivial. Alors $\operatorname{deg}(v_{\ast}) \geq 3$. Comme chaque composante connexe de $\overline{X}-\{v_{\ast}\}$ contient au moins une feuille, et qu'il existe un sommet de groupe associé non trivial et qui ne soit pas une feuille, $\operatorname{deg}(v_{\ast})=3$. De plus, puisqu'il y a  exactement $4$ sommets dans le graphe de groupes associés non triviaux, chaque composante connexe de $\overline{X}-\{v_{\ast}\}$ contient exactement une feuille. Donc $v_{\ast}$ est relié à exactement $2$ feuilles et $w$ est relié à une seule feuille et à $v_{\ast}$. Or le cardinal du groupe des automorphismes d'un tel graphe est égal à $2$. Comme le noyau du morphisme $G \to \Aut_{gr}(X)$ est de cardinal $4$, ceci contredit le fait que $G$ est isomorphe à $\mathfrak{S}_4$.

\bigskip

Supposons que $\overline{X}$ contienne un seul sommet, noté $w$, de groupe associé non trivial et qui ne soit pas une feuille. Si $v_{\ast}$ est une feuille, alors le graphe possède exactement $3$ feuilles, dont l'une est le point base. De ce fait, comme tout automorphisme de $\overline{X}$ est induit par son action sur les feuilles, le groupe des automorphismes d'un tel graphe pointé est de cardinal $2$. Comme le noyau du morphisme $G \to \Aut_{gr}(X)$ est de cardinal $4$, ceci contredit le fait que $G$ est isomorphe à $\mathfrak{S}_4$.

Supposons alors que le point base $v_{\ast}$ ne soit pas une feuille. Par les cas précédents, $v_{\ast}=w$. Comme le nombre de sommets de groupes non trivial est exactement $4$, et que tout sommet de groupe associé trivial est de degré au moins $3$, le graphe $\overline{X}$ contient au plus un sommet de groupe associé trivial. Le cas où le nombre de sommets de groupe associé trivial est égal à $1$ n'est pas possible car alors le cardinal du groupe des automorphismes d'un tel graphe est égal à $2$, contredisant le fait que le noyau du morphisme $G \to \Aut_{gr}(X)$ est de cardinal $4$ et que $G$ est isomorphe à $\mathfrak{S}_4$. 

Dans le cas où le nombre de sommets de groupe associé trivial est nul, on voit que $\mathcal{X}$ est une $F$-étoile. Or, par la remarque~\ref{rmq twists Aut}, le cardinal du noyau du morphisme $G \to \Aut_{gr}(X)$ est égal à $8$, d'où une contradiction.

\bigskip

En conclusion, le morphisme $G \to \Aut_{gr}(X)$ est également injectif dans le cas où $\mathcal{X}$ appartient à $\PA(W_4)$ et $n=4$. La suite de la démonstration est alors identique à la proposition~\ref{Sn fixe 0 étoile}.}
\hfill\qedsymbol

\end{rmq}

\bigskip

Nous démontrons à présent un résultat similaire pour les $F$-étoiles. Pour cela, nous avons besoin du lemme suivant.

\begin{lem}\label{maximalité k}
Soient $n \geq 4$ un entier et $\mathcal{X}$ un point de l'épine de $\PO(W_n)$. On note $X$ un représentant de $\mathcal{X}$ et $\overline{X}$ le graphe sous-jacent à $X$. Soit $k$ l'entier tel que le noyau du morphisme naturel $\Stab_{\Out(W_n)}(\mathcal{X}) \to \Aut_{gr}(X)$ soit isomorphe à $F^k$.
Alors $k \leq n-2$. Par ailleurs, $k=n-2$ si, et seulement si, l'ensemble $V\overline{X}$ des sommets de $\overline{X}$ est de cardinal $n$.
\end{lem}

\dem
Supposons que $|V\overline{X}|>n$. Soient $v$ un sommet de groupe associé trivial et $e$ une arête de $X$ reliant $v$ à un sommet $w$. Une telle arête existe car $\overline{X}$ est connexe et le nombre de sommets de $\overline{X}$ de groupe non trivial est égal à $n$.

\bigskip

\noindent{\bf Affirmation. } Soient $Y$ le graphe de groupes marqué obtenu à partir de $X$ en contractant l'arête $e$ et $\mathcal{Y}$ sa classe d'équivalence dans l'épine de $\PO(W_n)$. Alors le noyau du morphisme naturel $\Stab_{\Out(W_n)}(\mathcal{Y}) \to \Aut_{gr}(Y)$ est isomorphe à $F^l$, avec $l = k$ si le groupe associé à $w$ est trivial, et $l \geq k+1$ sinon.

\bigskip

\dem Si le groupe associé à $w$ est trivial, alors contracter l'arête $e$ ne modifie pas le degré des sommets dont le groupe associé est non trivial. Donc, dans ce cas, $k=l$. Supposons maintenant que le groupe associé à $w$ ne soit pas trivial. Notons $\widebar{vw}$ le sommet obtenu en contractant $e$. Le groupe associé à $\widebar{vw}$ est non trivial. Alors, puisque, par hypothèse, $\operatorname{deg}(v) \geq 3$, nous avons~:

\begin{equation*}
\operatorname{deg}(\widebar{vw}) =\operatorname{deg}(v) + \operatorname{deg}(w) -2 \geq \operatorname{deg}(w) +1.
\end{equation*}
Ainsi, dans ce cas, $l \geq k+1$. 
\hfill\qedsymbol

\medskip

De ce fait, si $|V\overline{X}|>n$, il existe une arête reliant un sommet de groupe associé trivial et un sommet de groupe associé non trivial. Par l'affirmation précédente, l'entier $k$ associé au morphisme $\Stab_{\Out(W_n)}(\mathcal{X}) \to \Aut_{gr}(X)$ n'est pas maximal. 

\bigskip

Ainsi, pour calculer la borne maximale de $k$, nous pouvons supposer que $\overline{X}$ possède $n$ sommets, tous de groupe associé non trivial. Donc, $$\sum_{v \in V\overline{X}} \operatorname{deg}(v)=2|E\overline{X}|=2n-2~,$$ la dernière égalité provenant du fait que $\overline{X}$ soit un arbre. Ainsi,

\begin{equation*}
k=\sum\limits_{v \in V\overline{X}}(\operatorname{deg}(v)-1)=\sum_{v \in V\overline{X}} \operatorname{deg}(v)-n=2n-2-n=n-2.
\end{equation*}

Donc, $k \leq n-2$, et si $|V\overline{X}|=n$, alors $k=n-2$. 

\medskip

Supposons maintenant que $k=n-2$. Par l'affirmation précédente, la procédure de contraction présentée fait croître strictement $k$ lorsque l'on contracte une arête reliant un sommet de groupe associé trivial et un sommet de groupe associé non trivial. Donc $\overline{X}$ ne peut pas contenir de sommets ayant un groupe associé trivial. Donc le cardinal de $V\overline{X}$ est égal à $n$.
\hfill\qedsymbol

\begin{rmq}\label{maximalité k Aut}
{\rm Dans le cas de $\PA(W_n)$, soit $\mathcal{X}$ un point de l'épine de $\PA(W_n)$. On note $X$ un représentant de $\mathcal{X}$ et $\overline{X}$ le graphe sous-jacent à $X$. Soit $k$ l'entier tel que le noyau du morphisme naturel $\Stab_{\Aut(W_n)}(\mathcal{X}) \to \Aut_{gr}(X)$ soit isomorphe à $F^k$. Alors une démonstration identique au lemme~\ref{maximalité k} montre que $k \leq n-1$ avec égalité si, et seulement si, $|V\overline{X}|=n$. }
\end{rmq}

\bigskip

Nous pouvons maintenant montrer le résultat suivant concernant les stabilisateurs de $F$-étoiles dans $\Out(W_n)$.

\begin{prop}\label{Wn fixe Z2 étoile}
\begin{enumerate}
\item Soit $n \geq 4$ un entier. Le cardinal maximal d'un sous-groupe fini de $\Out(W_n)$ est $2^{n-2}(n-1)!\;$.
\item Supposons $n \geq 5$. Soient $G$ un sous-groupe de $\Out(W_n)$, et $\mathcal{X}$ un point de l'épine de $\PO(W_n)$ fixé par $G$. On note $X$ un représentant de $\mathcal{X}$ et $\overline{X}$ le graphe sous-jacent à $X$. Si $\overline{X}$ possède $n$ feuilles, alors $|G|< 2^{n-2}(n-1)!\;$.
\item Supposons $n \geq 4$. Soient $G$ un sous-groupe de $\Out(W_n)$ isomorphe à $F^{n-2} \rtimes \mathfrak{S}_{n-1}$, et $\mathcal{X}$ un point de l'épine de $\PO(W_n)$ fixé par $G$. On note $X$ un représentant de $\mathcal{X}$ et $\overline{X}$ le graphe sous-jacent à $X$. Si le nombre de feuilles de $\overline{X}$ est $n-1$, alors $\mathcal{X}$ est une $F$-étoile.
\item Supposons $n \geq 5$. Soit $G$ un sous-groupe de $\Out(W_n)$ isomorphe à $F^{n-2} \rtimes \mathfrak{S}_{n-1}$. Alors $G$ est le stabilisateur d'une unique $F$-étoile.
\end{enumerate}
\end{prop}

\dem Si $\mathcal{X}$ est un élément de l'épine de $\PO(W_n)$, nous noterons $X$ un représentant de $\mathcal{X}$. Nous noterons également $\overline{X}$ le graphe sous-jacent à $X$ et $L$ l'ensemble des feuilles de $\overline{X}$. Puisque $\overline{X}$ est un arbre, tout automorphisme de $\overline{X}$ est entièrement déterminé par son action sur les feuilles. Donc le morphisme de restriction de $\Aut_{gr}(X)$ dans $\Bij(L)$ est injectif.

Montrons l'assertion~$(1)$. Puisque tout sous-groupe fini de $\Out(W_n)$ fixe un point de l'épine de $\PO(W_n)$ par la proposition~\ref{point fixe outre espace}, il suffit de montrer que, pour $\mathcal{X}$ un point de l'épine de l'outre-espace, $|\Stab_{\Out(W_n)}(\mathcal{X})|\leq 2^{n-2}(n-1)!\;$. D'après la proposition~\ref{prop Levitt}, il existe un entier $k$ tel que le noyau du morphisme naturel $\Stab_{\Out(W_n)}(\mathcal{X}) \to \Aut_{gr}(X)$ soit isomorphe à $F^k$.
De ce fait, $|\Stab_{\Out(W_n)}(\mathcal{X})|\leq 2^k|\Aut_{gr}(X)|$. 

Nous distinguons deux cas, selon le cardinal de $L$.
\begin{itemize}
\item Supposons que $|L| \leq n-1$. Alors $\Aut_{gr}(X)$, qui s'injecte dans $\Bij(L)$, s'injecte dans $\mathfrak{S}_{n-1}$. Ainsi,
\begin{equation*}
|\Stab_{\Out(W_n)}(\mathcal{X})|\leq 2^k|\Aut_{gr}(X)| \leq 2^k(n-1)! \leq 2^{n-2}(n-1)!\;,
\end{equation*}
où la dernière inégalité découle du lemme~\ref{maximalité k}.

\item Supposons que $|L|=n$. Alors tous les sommets ayant des groupes associés non triviaux sont des feuilles. Ainsi, $k=0$ par la proposition~\ref{prop Levitt}. Puisque $\Bij(L)$ est isomorphe à $\mathfrak{S}_n$, nous avons
\begin{equation*}
|\Stab_{\Out(W_n)}(\mathcal{X})|\leq |\Aut_{gr}(X)| \leq n!\;.
\end{equation*}
Or puisque $n \geq 4$, nous avons $n \leq 2^{n-2}$, donc $n! \leq 2^{n-2}(n-1)!$, ce qui conclut.
\end{itemize}

Donc, pour tout sous-groupe fini $G$ de $\Out(W_n)$, l'ordre de $G$ est au plus $2^{n-2}(n-1)!\;$. Cette borne est atteinte par le groupe $U_n=\left\langle [\tau_1],\ldots, [\tau_{n-2}], [\sigma_{1,n}] \right\rangle$. qui est isomorphe au produit semi-direct $F^{n-2} \rtimes \mathfrak{S}_{n-1}$.

\bigskip

Soient $n \geq 5$ et $G$, $\mathcal{X}$ et $X$ comme dans l'énoncé de l'assertion~$(2)$. Par la proposition~\ref{prop Levitt}, il existe un entier $k$ tel que le noyau du morphisme naturel $G \to \Aut_{gr}(X)$ soit isomorphe à $F^k \cap G$. Puisque $\overline{X}$ possède $n$ feuilles, par la proposition~\ref{prop Levitt}, l'entier $k$ est nul. De ce fait, le groupe $G$ s'injecte dans $\Aut_{gr}(X)$,  qui s'injecte dans $\Bij(L)$. Donc $|G| \leq n!\;$. Or $2^{n-2}(n-1)! \leq n!$ implique que $n \leq 4$. D'où $|G| < 2^{n-2}(n-1)!\;$.

\bigskip

Soient $n \geq 4$ et $G$, $\mathcal{X}$ et $X$ comme dans l'énoncé de $(3)$. Comme $G$ est de cardinal maximal parmi les sous-groupes finis de $\Out(W_n)$, nous avons $G=\Stab_{\Out(W_n)}(\mathcal{X})$. Donc, par la proposition~\ref{prop Levitt}, il existe un entier $k$ tel que le noyau du morphisme naturel $G \to \Aut_{gr}(X)$ soit isomorphe à $F^k$. Ainsi, puisque $\Aut_{gr}(X)$ s'injecte dans $\Bij(L)$ et que ce dernier est isomorphe à $\mathfrak{S}_{n-1}$, on voit que $|G| \leq 2^k(n-1)!\;$. Comme $k \leq n-2$ par le lemme~\ref{maximalité k}, et puisque $|G|=2^{n-2}(n-1)!$, on a nécessairement $k=n-2$. Le lemme~\ref{maximalité k} donne alors que $\overline{X}$ possède exactement $n$ sommets. De ce fait, $\overline{X}$ possède $n-1$ feuilles et $n$ sommets. Par ailleurs, on voit également que $\Aut_{gr}(X)$ est isomorphe à $\Bij(L)$. De ce fait, toutes les arêtes de $\overline{X}$ ont la même longueur. Donc $\mathcal{X}$ est une $F$-étoile.

\bigskip

Supposons enfin que $n \geq 5$ et que $G$ soit un sous-groupe de $\Out(W_n)$ isomorphe à $F^{n-2} \rtimes \mathfrak{S}_{n-1}$. Par la proposition~\ref{point fixe outre espace}, le groupe $G$ fixe un point $\mathcal{X}$ de l'épine de l'outre-espace. Comme $G$ est de cardinal maximal parmi les sous-groupes finis de $\Out(W_n)$, nous avons $G=\Stab_{\Out(W_n)}(\mathcal{X})$. Donc, par la proposition~\ref{prop Levitt}, il existe un entier $k$ tel que le noyau du morphisme naturel $G \to \Aut_{gr}(X)$ soit isomorphe à $F^k$.

\bigskip

\noindent{\bf Affirmation. } L'arbre $\overline{X}$ possède exactement $n-1$ feuilles.

\bigskip

\dem L'assertion~$(2)$ dit que $\overline{X}$ possède au plus $n-1$ feuilles. Nous avons $$|G|=2^{n-2}(n-1)! \leq 2^{k}|\Aut_{gr}(X)|\leq 2^{n-2}|\Aut_{gr}(X)|~;$$ où la dernière égalité provient du lemme~\ref{maximalité k}. Donc $|\Aut_{gr}(X)| \geq (n-1)!\;$. Ainsi, puisque $\overline{X}$ possède au plus $n-1$ feuilles, le groupe $\Bij(L)$, dans lequel s'injecte $\Aut_{gr}(X)$, est isomorphe à $\mathfrak{S}_{n-1}$. Donc le cardinal de $L$ est $n-1$. \hfill\qedsymbol

\bigskip

De ce fait, $\overline{X}$ possède $n-1$ feuilles. Par l'assertion~$(3)$, $\mathcal{X}$ est une $F$-étoile dans l'épine de $\PO(W_n)$. Par le corollaire~\ref{points fixes connexes}, l'ensemble des points fixes de $G$ est connexe. Puisque l'ensemble des $F$-étoiles est discret dans l'épine de $\PO(W_n)$, on conclut que $G$ fixe une unique $F$-étoile dans l'épine de $\PO(W_n)$.
\hfill\qedsymbol

\begin{rmq}\label{Wn fixe Z2 étoile Aut}
{\rm Dans le cas de $\Aut(W_n)$, soient $G$ un sous groupe fini de $\Aut(W_n)$ et $\mathcal{X}$ un point de l'épine de $\PA(W_n)$ fixé par $G$. On note $X$ un représentant de $\mathcal{X}$ et $\overline{X}$ le graphe pointé sous-jacent à $X$. 

\smallskip

{\it \noindent{$(1)$ } Si $n \geq 4$, le cardinal de $G$ est plus petit que $2^{n-1}(n-1)!\,.$ }

La démonstration pour le cas où le nombre de feuilles de $\overline{X}$ est plus petit que $n-1$ est identique à celle de la proposition~\ref{Wn fixe Z2 étoile}~$(1)$ en utilisant cette fois la remarque~\ref{maximalité k Aut}. Dans le cas où le nombre de feuilles est égal à $n$, le noyau du morphisme naturel $G \to \Aut_{gr}(X)$ est de cardinal plus petit que $2$ par la remarque~\ref{rmq twists Aut}, donc $|G| \leq 2n! \leq 2^{n-1}(n-1)!$ car $n \geq 4$.

\smallskip

{\it \noindent{$(2)$}  Si $n \geq 5$ et si $\overline{X}$ possède $n$ feuilles, alors $|G|< 2^{n-1}(n-1)!\;$. }

En effet, par la remarque~\ref{rmq twists Aut}, le cardinal du noyau du morphisme $G \to \Aut_{gr}(X)$ est plus petit que $2$, donc $|G| \leq 2n! < 2^{n-1}(n-1)!$ car $n \geq 5$.

\smallskip

{\it \noindent{$(3)$ } Si $n \geq 4$, si $G$ est isomorphe à $F^{n-1} \rtimes \mathfrak{S}_{n-1}$ et si $\overline{X}$ possède au plus $n-1$ feuilles, alors $\mathcal{X}$ est une $F$-étoile. }

En effet, une démonstration identique à celle de la proposition~\ref{Wn fixe Z2 étoile}~$(3)$ montre que $\overline{X}$ possède $n-1$ feuilles et $n$ sommets. Montrons alors que le point base est le centre de $X$. Ceci découle du fait que le groupe des automorphismes de $\overline{X}$ est isomorphe à $\mathfrak{S}_{n-1}$ car le noyau du morphisme $G \to \Aut_{gr}(X)$ est isomorphe à $F^{n-1}$ et que $G$ est isomorphe à $F^{n-1} \rtimes \mathfrak{S}_{n-1}$.

\smallskip

{\it \noindent{$(4)$ } Si $n \geq 4$ et si $G$ est isomorphe à $F^{n-1} \rtimes \mathfrak{S}_{n-1}$, tout point de l'épine de $\PA(W_n)$ fixé par $G$ est une $F$-étoile.}

En effet, l'existence d'une $F$-étoile fixée par $G$ lorsque $n \geq 5$ se déduit des faits précédents.

Dans le cas où $n=4$, soit $\mathcal{X}$ un point de l'épine de l'outre-espace fixé par $G$. Soient $X$ un représentant de $\mathcal{X}$ et $\overline{X}$ le graphe sous-jacent à $X$. On note $L$ l'ensemble des feuilles de $\overline{X}$. Si $\overline{X}$ possède au plus $n-1$ feuilles, alors, par le fait précédent, $\mathcal{X}$ est une $F$-étoile. Supposons que $\overline{X}$ possède exactement $n$ feuilles. Alors la remarque~\ref{rmq twists Aut} montre que le noyau du morphisme naturel $G \to \Aut_{gr}(X)$ est de cardinal au plus $2$. Il ne peut pas être injectif car le cardinal de $G$ est égal à $48$ alors que le groupe $\Aut_{gr}(X)$ s'injecte dans $\Bij(L)$ de cardinal égal à $24$. Donc le noyau du morphisme $G \to \Aut_{gr}(X)$ est de cardinal égal à $2$. Ainsi, le point base de $\overline{X}$ est une feuille. Or, puisque $\Aut_{gr}(X)$ s'injecte dans $\Bij(L)$ et que l'image du morphisme $G \to \Aut_{gr}(X)$ est de cardinal égal à $24$, on voit que $\Aut_{gr}(X)$ est isomorphe à $\Bij(L)$. Ceci contredit le fait que le point base de $\overline{X}$ est une feuille. En conclusion, $\overline{X}$ possède au plus $n-1$ feuilles. Donc $\mathcal{X}$ est une $F$-étoile. La démonstration de l'unicité de la $F$-étoile fixée par $G$ est alors identique à celle de la démonstration de la proposition~\ref{Wn fixe Z2 étoile}~$(4)$.}
\end{rmq}

\begin{lem}\label{maximalité de Sn-1}
Soit $n$ un entier.
\begin{enumerate}
\item Supposons que $n \geq 5$. Soit $G$ un sous-groupe de $\mathfrak{S}_n$ isomorphe à $\mathfrak{S}_{n-1}$. Il existe un automorphisme de $\mathfrak{S}_n$ envoyant $G$ sur $\{f \in \Bij(\{1,\ldots,n\})\; : \; f(n)=n\}$.
\item Si $n \geq 4$ et $n \neq 6$ et si $G$ est un sous-groupe de $\Bij(\{1,\ldots,n\})$ isomorphe à $\mathfrak{S}_{n-1}$, alors il existe un entier $i \in \{1,\ldots,n\}$ tel que $G=\{f \in \Bij(\{1,\ldots,n\})\; : \; f(i)=i\}$.
\end{enumerate}

\end{lem}

\dem \noindent{(1) } L'action de $\mathfrak{S}_n$ sur $\mathfrak{S}_n/G$ par multiplication à gauche est un morphisme de groupes $\phi \colon \mathfrak{S}_n \to \mathrm{Bij}(\mathfrak{S}_n/G)$. Le noyau de ce morphisme est un sous-groupe distingué de $\mathfrak{S}_n$ inclus dans $G$. Or, $G$ est d'indice $n$. Donc, étant donné que $n \geq 5$, le noyau de ce morphisme est trivial. Donc, puisque les groupes $\mathfrak{S}_n$ et $\mathrm{Bij}(\mathfrak{S}_n/G)$ ont même cardinal fini, le morphisme $\phi$ est un isomorphisme. Soit $\widetilde{\psi} \colon \mathfrak{S}_n/G \to \{1,\ldots,n\}$ une bijection envoyant $\{G\}$ sur $n$, et $\psi \colon \mathrm{Bij}(\mathfrak{S}_n/G) \to \mathfrak{S}_n$ l'isomorphisme induit par $\widetilde{\psi}$. Alors $\psi \circ \phi$ est un automorphisme de $\mathfrak{S}_n$ envoyant $G$ sur le sous-groupe de $\mathfrak{S}_n$ fixant $n$.

\medskip

\noindent{(2) } Nous commençons par traiter le cas où $n=4$. Il découle d'une inspection des sous-groupes de $\mathfrak{S}_4$ isomorphes à $\mathfrak{S}_3$. En effet, $\mathfrak{S}_4$ possède exactement $4$ sous-groupes isomorphes à $\mathfrak{S}_3$. Donc, il existe un entier $i \in \{1,2,3,4\}$ tel que $G=\{f \in \Bij(\{1,\ldots,n\})\; : \;f(i)=i\}$.

Supposons maintenant que $n \geq 5$ et que $n \neq 6$. Par le premier point du lemme, il existe un automorphisme $\phi$ de $\mathfrak{S}_n$ envoyant $G$ sur $\{f \in \Bij(\{1,\ldots,n\})\; : \; f(n)=n\}$. Or, si $n \neq 6$, tout automorphisme de $\mathfrak{S}_n$ est intérieur. Comme les automorphismes intérieurs préservent le fait d'être le stabilisateur d'un entier, il existe un entier $i \in \{1,\ldots,n\}$ tel que $G=\{f \in \Bij(\{1,\ldots,n\})\; : \; f(i)=i\}$.
\hfill\qedsymbol

\bigskip

Étudions les points fixes du groupe $B_n$ dans l'épine de l'outre-espace de $W_n$.

\begin{prop}\label{unicité Z2 étoile}
Soient $n \geq 4$ et $B_n=\left\langle [\tau_1],\ldots, [\tau_{n-2}] \right\rangle$. 
\begin{enumerate}
\item Les seuls sommets fixés par $B_n$ dans l'épine de l'outre-espace de $W_n$ sont des $\{0\}$-étoiles et des $F$-étoiles.
\item Le groupe $B_n$ fixe une unique $F$-étoile et une unique $\{0\}$-étoile.
\end{enumerate}
\end{prop}

\noindent{\bf Remarque. }{\rm La proposition~\ref{unicité Z2 étoile} diffère des propositions~\ref{Sn fixe 0 étoile} et \ref{Wn fixe Z2 étoile} car elle porte uniquement sur un sous-groupe \emph{particulier} de $\Out(W_n)$. Nous ne savons pas si le résultat reste vrai pour un sous-groupe de $\Out(W_n)$ isomorphe à $\mathfrak{S}_{n-1}$ quelconque.}

\bigskip

\dem \noindent{(1) } Soient $\mathcal{X}$ un sommet de l'épine de $\PO(W_n)$ fixé par $B_n$ et $X$ un représentant de $\mathcal{X}$. Soient $\overline{X}$ le graphe sous-jacent à $X$, $L$ l'ensemble des feuilles de $\overline{X}$ et $v_1,\ldots,v_n$ les sommets de $\overline{X}$ dont les groupes associés sont non triviaux. Par la proposition~\ref{prop Levitt}, il existe un entier $k$ tel que le noyau du morphisme naturel \mbox{$B_n \to \Aut_{gr}(X)$} soit isomorphe à $F^k \cap B_n$.
Or, ce noyau est un sous-groupe de $F^k$, et ce dernier est engendré par des twists. Pour tout $i \in \{1,\ldots,n\}$, soit $y_i$ l'antécédent par le marquage de $X$ du générateur du groupe associé à $v_i$. Pour tout $i \in \{1,\ldots,n\}$, les compositions de twists contenues dans $F^k \cap B_n$ préservent la classe de conjugaison dans $W_n$ de $y_i$ alors que les permutations du groupe engendré par $\{[\tau_1],\ldots, [\tau_{n-2}]\}$ ne préservent pas ces dernières. De ce fait, nous avons $F^k \cap B_n=\{1\}$. 

Le groupe $\Aut_{gr}(X)$ s'injecte dans $\mathrm{Bij}(L)$. Par ailleurs, étant donné que le morphisme $\phi \colon B_n \to \Aut_{gr}(X)$ est injectif, et que $B_n$ est isomorphe à $\mathfrak{S}_{n-1}$, nous avons $|L| \geq n-1$. De plus, chaque feuille ayant un groupe associé non trivial, nous avons $|L| \leq n$. Donc $L \in \{n-1,n\}$. Examinons les deux cas possibles.

Si $|L|=n-1$, alors $\Aut_{gr}(X)$ est isomorphe à $\mathrm{Bij}(L)$. Montrons que $\mathcal{X}$ est une $F$-étoile. Soit $v$ un sommet qui n'est pas une feuille à distance maximale du centre de $\overline{X}$. L'hypothèse de maximalité sur $v$ implique qu'il y a au plus un sommet adjacent à $v$ qui n'est pas une feuille, car sinon nous pourrions trouver un sommet $w$ adjacent à $v$ qui ne serait pas une feuille et qui serait à distance strictement plus grande du centre que $v$. De ce fait, $v$ est adjacent à au moins $\operatorname{deg}(v)-1$ feuilles. 

Si le groupe associé à $v$ est non trivial, alors $v$ est fixé par $B_n$ car c'est le seul sommet de $\overline{X}$ qui soit de groupe associé non trivial et qui ne soit pas une feuille. Donc puisque $B_n$ est isomorphe à $\Aut_{gr}(X)$, le sommet $v$ est fixé par $\Aut_{gr}(X)$. Enfin, puisque tout élément de $\Bij(L)$ est induit par un élément de $\Aut_{gr}(X)$, le sommet $v$ est adjacent à toutes les feuilles et $\mathcal{X}$ est une $F$-étoile. 

Si $v$ est un sommet de groupe trivial, alors, par hypothèse, $\operatorname{deg}(v) \geq 3$. De ce fait, $v$ est adjacent à au moins deux feuilles, notées $v_1$ et $v_2$. Soit $w$ une feuille de $\overline{X}$ distincte de $v_1$ et $v_2$. Puisqu'il existe un automorphisme de $\overline{X}$ envoyant $v_1$ sur $w$ et fixant $v_2$, alors, nécessairement, $w$ est adjacent à $v$. Donc $v$ est adjacent à toutes les feuilles. Ceci n'est pas possible car alors $X$ contiendrait uniquement $n-1$ sommets de groupe associé non trivial. Donc $v$ est nécessairement un sommet de groupe associé non trivial et $\mathcal{X}$ est une $F$-étoile.

\bigskip

Supposons que $|L|=n$. Montrons alors que $\mathcal{X}$ est une $\{0\}$-étoile. Le groupe $\Aut_{gr}(X)$ s'injecte dans $\Bij(L)$ qui est isomorphe à $\mathfrak{S}_n$. Par ailleurs, puisque $B_n$ s'injecte dans $\Aut_{gr}(X)$, l'image de $\Aut_{gr}(X)$ dans $\Bij(L)$ contient un sous-groupe de $\Bij(L)$ isomorphe à $\mathfrak{S}_{n-1}$.

Soit $H$ l'image de $B_n$ dans $\Aut_{gr}(X)$. Par le lemme~\ref{maximalité de Sn-1}~$(2)$, si $n \neq 6$, il existe une feuille $v_1$ de $\overline{X}$ telle que l'image de $H$ dans $\Bij(L)$ soit égale à $\Stab_{\Bij(L)}(v_1)$. Soit $v$ le sommet adjacent à $v_1$. Puisque $v$ n'est pas une feuille, $\operatorname{deg}(v) \geq 3$. Ou bien $v$ est adjacent à une autre feuille distincte de $v_1$, ou bien $v$ est adjacent à une unique feuille. 

Si $v$ est adjacent à une unique feuille, il existe dans $\overline{X}$ des feuilles de $L-\{v_1\}$ à distance au moins $4$. Soient $w_1$ et $w_2$ deux telles feuilles distinctes de $v_1$, telles que $w_1$ soit à distance maximale du centre et que $w_2$ soit une feuille distincte de $v_1$ à distance maximale de $w_1$. Puisque la valence de tout sommet de groupe associé trivial est au moins $3$, il existe une feuille $w_3$ à distance $2$ de $w_2$. Or l'image de $H$ dans $\Bij(L)$ est égale à $\Stab_{\Bij(L)}(v_1)$. Donc il existe un automorphisme de $\overline{X}$ fixant $w_3$ et envoyant $w_2$ sur $w_1$, ce qui n'est pas possible par hypothèse sur $w_1$ et $w_2$. 

Donc $v$ est adjacent à une feuille distincte de $v_1$, que l'on note $v_2$. Soit $w$ une feuille de $\overline{X}$ distincte de $v_1$ et $v_2$. Étant donne qu'il existe un automorphisme de $\overline{X}$ envoyant $v_2$ sur $w$ et fixant $v_1$, le sommet $w$ est à distance $2$ de $v_2$. En particulier, $\mathcal{X}$ est une $\{0\}$-étoile.

\bigskip

Traitons maintenant le cas où $n=6$. On numérote de $1$ à $6$ les feuilles. Une construction explicite d'un représentant de l'unique automorphisme extérieur non trivial de $\mathfrak{S}_6$ (cf. \cite{Miller58}) donne que l'unique (à conjugaison près) sous-groupe de $\Bij(L)$ isomorphe à $\mathfrak{S}_5$ et qui
ne soit pas un stabilisateur de feuille est le groupe $$H=\left\langle (1 \,2)(3 \,4)(5 \, 6),(1 \,6)(2 \,4)(3 \, 5),(1 \,4)(2 \,3)(5 \, 6),(1 \,6)(2 \,5)(3 \, 4)\right\rangle.$$ Supposons alors que $H$ soit inclus dans l'image de $\Aut_{gr}(X)$ dans $\Bij(L)$. Le groupe $H$ agit transitivement sur les feuilles de $\overline{X}$. De ce fait, tous les sommets reliés à des feuilles sont adjacents à un même nombre $k$ de feuilles. Les seules valeurs possibles pour $k$ sont $k \in \{1,2,3,6\}$.  Le cas où $k=1$ n'est pas possible car tout sommet qui n'est pas une feuille est de degré au moins $3$ (tous les sommets dont les groupes associés sont non triviaux sont des feuilles). De plus, $k \neq 3$ car le groupe des automorphismes d'un tel graphe ne pourrait contenir simultanément les permutations $(1 \,2)(3 \,4)(5 \, 6)$, $(1 \,6)(2 \,4)(3 \, 5)$ et $(1 \,4)(2 \,3)(5 \, 6)$. Enfin, $k \neq 2$ car alors $\overline{X}$ posséderait $3$ sommets adjacents à $2$ feuilles. Cependant le groupe des automorphismes d'un tel graphe ne pourrait contenir simultanément les permutations $(1 \,2)(3 \,4)(5 \, 6)$, $(1 \,6)(2 \,4)(3 \, 5)$ et $(1 \,6)(2 \,5)(3 \, 4)$. Donc $k=6$ et $X$ est une $\{0\}$-étoile.

Ainsi, $B_n$ fixe uniquement des $\{0\}$-étoiles et des $F$-étoiles. 

\medskip

\noindent{(2) } Montrons maintenant que $B_n$ fixe une unique $F$-étoile. Soit $X$ le graphe de groupes marqué dont le graphe sous-jacent possède $n$ sommets, notés $v_1,\ldots,v_n$, tel que les feuilles du graphe sous-jacent soient $v_1,\ldots,v_{n-1}$, et tel que pour tout $i \in \{1,\ldots,n\}$, l'image réciproque par le marquage du générateur du groupe associé à $v_i$ soit $x_i$. Soit $\mathcal{X}$ la classe d'équivalence de $X$. Alors $\mathcal{X}$ est une $F$-étoile et le stabilisateur de $\mathcal{X}$ est $U_n$. Puisque $B_n \subseteq U_n$, ceci montre l'existence.

Montrons maintenant l'unicité. Soit $\mathcal{Y}$ une autre $F$-étoile fixée par $B_n$. On note $Y$ un représentant de $\mathcal{Y}$. Par le corollaire~\ref{points fixes connexes}, il existe dans $\mathrm{Fix}_{K_n}(B_n)$ un chemin continu de $\mathcal{X}$ vers $\mathcal{Y}$. Puisque deux $F$-étoiles distinctes ne sont pas reliées par une arête dans l'épine de $\PO(W_n)$, et puisque tout sommet de $\mathrm{Fix}_{K_n}(B_n)$ est une $\{0\}$-étoile ou une $F$-étoile, ce chemin passe par une $\{0\}$-étoile adjacente à $\mathcal{X}$.

\bigskip

\noindent{\bf Affirmation. } Soient $\mathcal{Z}$ une $\{0\}$-étoile adjacente à $\mathcal{X}$ et $Z$ un représentant de $\mathcal{Z}$. On note $\overline{Z}$ le graphe sous-jacent à $Z$ et $v_1,\ldots,v_n$ les sommets de $\overline{Z}$ dont les groupes associés sont non triviaux. Alors l'image réciproque par le marquage de $Z$ des générateurs des groupes associés aux sommets $v_1,\ldots,v_n$ est, à conjugaison près~: $$\{x_n^{\alpha_1}x_1x_n^{\alpha_1},\ldots,x_n^{\alpha_{n-1}}x_{n-1}x_n^{\alpha_{n-1}},x_n\}~,$$ avec $\alpha_i \in \{0,1\}$ pour tout $i \in \{1,\ldots,n-1\}$.

\bigskip

\dem Pour tout $i \in \{1,\ldots,n\}$, soit $y_i$ le générateur du groupe associé à $v_i$. Puisque $\mathcal{Z}$ est adjacente à $\mathcal{X}$, il existe une arête $e$ de $\overline{Z}$ telle que le graphe de groupes marqué $Z'$ dont le graphe sous-jacent est obtenu à partir de $\overline{Z}$ en contractant $e$ soit dans la classe $\mathcal{X}$. Quitte à renuméroter, on peut supposer que l'un des sommets de $e$ est $v_n$. Soient $T_X$ et $T_{Z'}$ les arbres de Bass-Serre associés à $X$ et $Z'$. Les graphes de groupes $X$ et $Z'$ étant équivalents, il existe un homéomorphisme $W_n$-équivariant $f \colon T_X \to T_{Z'}$. Soit $v$ le sommet de $T_X$ de stabilisateur $\left\langle x_n \right\rangle$. Alors $f(v)$ a pour stabilisateur $\left\langle x_n \right\rangle$. Par ailleurs, étant donné que les sommets adjacents à $v$ ont pour stabilisateurs $\left\langle x_1 \right\rangle, \ldots, \left\langle x_{n-1} \right\rangle,\left\langle x_nx_1x_n \right\rangle,\ldots, \left\langle x_nx_{n-1}x_n \right\rangle$, les sommets
adjacents à $f(v)$ ont pour stabilisateurs $\left\langle x_1 \right\rangle, \ldots, \left\langle x_{n-1} \right\rangle,\left\langle x_nx_1x_n \right\rangle,\ldots, \left\langle x_nx_{n-1}x_n \right\rangle$. Donc, tout sous-graphe fini et connexe de $T_{Z'}$ ayant $n$ sommets et $n-1$ feuilles et de centre $f(v)$ est tel que les stabilisateurs des feuilles sont $$\left\langle x_n^{\alpha_1}x_1x_n^{\alpha_1}\right\rangle,\ldots,\left\langle x_n^{\alpha_{n-1}}x_{n-1}x_n^{\alpha_{n-1}}\right\rangle,$$ avec  $\alpha_i \in \{0,1\}$ pour tout $i \in \{1,\ldots,n-1\}$. Ainsi, l'image réciproque par le marquage de $Z$ des générateurs des groupes associés aux sommets $v_1,\ldots,v_n$ est, à conjugaison près~: $$\left\langle x_n^{\alpha_1}x_1x_n^{\alpha_1}\right\rangle,\ldots,\left\langle x_n^{\alpha_{n-1}}x_{n-1}x_n^{\alpha_{n-1}}\right\rangle,$$ avec  $\alpha_i \in \{0,1\}$ pour tout $i \in \{1,\ldots,n-1\}$.
\hfill\qedsymbol
\bigskip

Ainsi, au vu de la description des $\{0\}$-étoiles adjacentes à $\mathcal{X}$, le groupe $B_n$ fixe une unique $\{0\}$-étoile adjacente à $\mathcal{X}$~: la $\{0\}$-étoile $Z$ telle que les antécédents par le marquage des générateurs des groupes de sommets non triviaux soient, à conjugaison près, $x_1,\ldots,x_n$. On note $\mathcal{Z}$ la classe d'équivalence de $Z$ et $\overline{Z}$ le graphe sous-jacent à $Z$.

\bigskip

Soit $\mathcal{Y}'$ une $F$-étoile adjacente à $\mathcal{Z}$. Notons $Y'$ un représentant de $\mathcal{Y}'$ et $\overline{Y}'$ le graphe sous-jacent à $Y$. Il existe une arête $e$ de $\overline{Z}$ telle que le graphe de groupes $Z'$ obtenu en contractant $e$ soit dans $\mathcal{Y}'$. Les antécédents par le marquage de $Y'$ des générateurs des groupes de sommets sont donc, à conjugaison près, $x_1,\ldots,x_n$. 

\bigskip

Ainsi, puisque $B_n$ permute les sommets de tout point de l'épine de $\PO(W_n)$ dont l'image réciproque par le marquage des groupes associés sont $\left\langle x_1\right\rangle,\ldots,\left\langle x_{n-1} \right\rangle$, on voit que l'unique $F$-étoile adjacente à $\mathcal{Z}$ fixée par $B_n$ est $\mathcal{X}$. Donc, $B_n$ fixe une unique $F$-étoile dans l'épine de $\PO(W_n)$.

\bigskip

Montrons enfin que $B_n$ fixe une unique $\{0\}$-étoile. Soit $Z$ le graphe de groupes marqué dont le graphe sous-jacent possède $n+1$ sommets, $n$ feuilles, notées $w_1,\ldots,w_n$, et tel que pour tout $i \in \{1,\ldots,n\}$, l'image réciproque par le marquage du générateur du groupe associé à $w_i$ soit $x_i$. Soit $\mathcal{Z}$ la classe d'équivalence de $Z$. Alors $\mathcal{Z}$ est une $\{0\}$-étoile et le stabilisateur de $\mathcal{Z}$ est $A_n$. Puisque $B_n \subseteq A_n$, ceci montre l'existence.
 
Montrons l'unicité. Soit $\mathcal{Y}$ une autre $\{0\}$-étoile fixée par $B_n$. Par le corollaire~\ref{points fixes connexes}, il existe un chemin continu dans $\mathrm{Fix}_{K_n}(B_n)$ de $\mathcal{Z}$ vers $\mathcal{Y}$. Au vu de l'assertion $(1)$ de la proposition, ce chemin passe uniquement par des $\{0\}$-étoiles et des $F$-étoiles. Or, $B_n$ fixe une unique $F$-étoile $\mathcal{X}$, et par la dernière affirmation, l'unique $\{0\}$-étoile adjacente à $\mathcal{X}$ et fixée par $B_n$ est $\mathcal{Z}$. Donc $B_n$ fixe une unique $\{0\}$-étoile dans l'épine de $\PO(W_n)$.

\hfill\qedsymbol

\begin{rmq}\label{Unicité Z2 étoile Aut}
{\rm Soit $n \geq 4$. Dans le cas de $\Aut(W_n)$, soit $\widetilde{B}_n=\left\langle \tau_1,\ldots, \tau_{n-2} \right\rangle$, qui est encore isomorphe à $\mathfrak{S}_{n-1}$. Soit $\mathcal{X}$ un point de l'épine de $\PA(W_n)$ fixé par $\widetilde{B}_n$. On note $X$ un représentant de $\mathcal{X}$ et $\overline{X}$ le graphe sous-jacent à $X$. 

\smallskip

{\it \noindent{$(1)$ } Soit $\mathcal{X}$ est une $F$-étoile, soit $\overline{X}$ possède $n$ feuilles et $n+1$ sommets. }

En effet, une démonstration identique à celle de la proposition~\ref{unicité Z2 étoile}~$(1)$ montre que le morphisme $\widetilde{B}_n \to \Aut_{gr}(X)$ est injectif, et que le nombre de feuilles de $\overline{X}$ est soit égal à $n-1$, soit égal à $n$. S'il est égal à $n-1$, une démonstration identique à celle de la proposition~\ref{unicité Z2 étoile}~$(1)$ montre que $\overline{X}$ possède $n$ sommets et $n-1$ feuilles. Comme le groupe $\Aut_{gr}(X)$ contient un sous-groupe isomorphe à $\mathfrak{S}_{n-1}$ et que $\overline{X}$ possède $n-1$ feuilles, on voit que, nécessairement, le point base de $X$ est son centre. Donc $\mathcal{X}$ est une $F$-étoile. 
Si le nombre de feuilles de $\overline{X}$ est égal à $n$, une démonstration identique à celle de la proposition~\ref{unicité Z2 étoile}~$(1)$ montre que $\overline{X}$ possède $n+1$ sommets et $n$ feuilles.

\smallskip

{\it \noindent{$(2)$ } Le groupe $\widetilde{B}_n$ fixe une unique $F$-étoile.}

En effet, il fixe une $F$-étoile car $\widetilde{B}_n$ est un sous-groupe de $\widetilde{U}_n=\left\langle \tau_1,\ldots, \tau_{n-2},\sigma_{1,n} \right\rangle$ et ce dernier est isomorphe à $F^{n-1} \rtimes \mathfrak{S}_{n-1}$. De ce fait, la remarque~\ref{Wn fixe Z2 étoile Aut}~$(4)$ permet de conclure. Nous appellerons $\mathcal{X}$ l'unique $F$-étoile fixée  par $\widetilde{U}_n$.

Pour l'unicité, soit $\mathcal{Y}$ une autre $F$-étoile fixée par $\widetilde{B}_n$. Puisque l'ensemble des $F$-étoiles dans l'épine de $\PA(W_n)$ n'est pas connexe, tout chemin continu entre $\mathcal{X}$ et $\mathcal{Y}$ et contenu dans l'ensemble des points fixes de $\widetilde{B}_n$ pour l'action de $\Aut(W_n)$ sur l'épine de $\PA(W_n)$ passe par un point $\mathcal{Z}$ ayant un représentant $Z$ de graphe sous-jacent possédant $n$ feuilles et $n+1$ sommets. Soit $\overline{Z}$ le graphe sous-jacent à $Z$, et $v_1,\ldots,v_n$ les feuilles de $\overline{Z}$. Une démonstration identique à celle de la première affirmation de la démonstration de la proposition~\ref{unicité Z2 étoile}~$(2)$ montre que l'image réciproque par le marquage de $Z$ des générateurs des groupes associés aux sommets $v_1,\ldots,v_n$ est respectivement ou bien $x_1,\ldots,x_{n-1},x_n$ ou bien $x_nx_1x_n,\ldots,x_nx_{n-1}x_n,x_n$. De plus, la description de $\widetilde{B}_n$ montre que le point base de $Z$ est contenu dans l'arête reliant le centre de $\overline{Z}$ et $v_n$. 

Soit maintenant $\mathcal{Z}'$ un sommet de l'épine de $\PA(W_n)$ fixé par $\widetilde{B}_n$, adjacent à $\mathcal{Z}$ et qui n'est pas une $F$-étoile. Puisque $Z'$ possède $n$ feuilles et $n+1$ sommets par le premier point de la remarque, un représentant $Z'$ de $\mathcal{Z}'$ est obtenu à partir de $Z$ en déplaçant le point base dans l'arête reliant le centre de $\overline{Z}$ et $v_n$. De ce fait, l'image réciproque par le marquage des générateurs des groupes associés aux feuilles de $\overline{Z}'$ sont les mêmes que pour $\mathcal{Z}$. 

Donc, pour conclure sur l'unicité de la $F$-étoile fixée par $\widetilde{B}_n$, il suffit d'étudier les $F$-étoiles fixées par $\widetilde{B}_n$ est adjacente à $\mathcal{Z}$. Soit $\mathcal{Y}'$ une $F$-étoile adjacente à $\mathcal{Z}$. Notons $Y'$ un représentant de $\mathcal{Y}'$ et $\overline{Y}'$ le graphe sous-jacent à $Y$. Il existe une arête $e$ de $\overline{Z}$ telle que le graphe de groupes $Z'$ obtenu en contractant $e$ soit dans $\mathcal{Y}'$. Les antécédents par le marquage de $Y'$ des générateurs des groupes de sommets sont donc, à conjugaison près, $x_1,\ldots,x_n$. Ainsi, puisque $\widetilde{B}_n$ permute transitivement les sommets de tout point de l'épine de $\PA(W_n)$ dont l'image réciproque par le marquage des groupes associés sont $\left\langle x_1\right\rangle,\ldots,\left\langle x_{n-1} \right\rangle$, on voit que l'unique $F$-étoile adjacente à $\mathcal{Z}$ fixée par $\widetilde{B}_n$ est $\mathcal{X}$. Donc, $\widetilde{B}_n$ fixe une unique $F$-étoile dans l'épine de $\PA(W_n)$. }
\end{rmq}

\section{Rigidité des automorphismes extérieurs d'un groupe de Coxeter universel}

Le but de cette partie est de démontrer le théorème~\ref{Théorème out coxeter}. Nous distinguons différents cas, selon la valeur de $n$. Soit $\alpha \in \Aut(\Out(W_n))$. 

\subsection{Démonstration dans le cas $n \geq 5$ et $n \neq 6$}\label{démonstration n=5}

Soit $\mathcal{X}_1$ la $\{0\}$-étoile fixée par le sous-groupe fini $A_n$ de $\Out(W_n)$ (l'unicité provient de la proposition~\ref{Sn fixe 0 étoile}). Alors, d'après la proposition~\ref{Sn fixe 0 étoile}, $\alpha(A_n)$ est le stabilisateur d'une unique $\{0\}$-étoile $\mathcal{X}_2$. Or $\Out(W_n)$ agit transitivement sur l'ensemble des $\{0\}$-étoiles, donc il existe $\psi \in \Out(W_n)$ tel que $\psi(\mathcal{X}_1)=\mathcal{X}_2$. Posons $\alpha_0=ad(\psi) \circ \alpha$, alors $\alpha_0(A_n)=ad(\psi) \circ \alpha(A_n)=A_n$. 

Puisque $\alpha_0|_{A_n}$ est un automorphisme de $A_n$, que $A_n$ est isomorphe à $\mathfrak{S}_n$ et que, pour $n \neq 6$, le groupe $\Out(\mathfrak{S}_n)$ est trivial, quitte à changer $\alpha_0$ dans sa classe d'automorphisme extérieur, on peut supposer que $\alpha_0|_{A_n}=\mathrm{id}_{A_n}$.

Maintenant, étant donné que $B_n \subseteq U_n$, nous avons $\alpha_0(B_n)=B_n \subseteq \alpha_0(U_n)$. Or par la proposition~\ref{unicité Z2 étoile}~$(2)$, $B_n$ fixe une unique $F$-étoile. Par ailleurs, le stabilisateur de cette $F$-étoile est $U_n$. Donc, puisque $\alpha_0(U_n)$ est également le stabilisateur d'une unique $F$-étoile par la proposition~\ref{Wn fixe Z2 étoile}~$(4)$, on obtient que $\alpha_0(U_n)=U_n$. 

Or $U_n$ est isomorphe au produit semi-direct $F^{n-2} \rtimes B_n$, et $B_n$ agit sur $F^{n-2}$ (vu comme le quotient de $F^{n-1}$ par son sous-groupe diagonal $F$) par permutation des facteurs. Soit $\sigma \in B_n$. On note $\mathrm{fix}(\sigma)$ l'ensemble des points fixes de $\sigma$ agissant par conjugaison dans $F^{n-2}$. Puisque, pour tout $\sigma \in B_n$, $\alpha_0(\sigma)=\sigma$, on voit que, pour tout $\sigma \in \{0\} \rtimes B_n$ et pour tout $g \in F^{n-2} \rtimes \{1\}$, $\alpha_0(\sigma g \sigma^{-1})=\sigma\alpha_{0}(g)\sigma^{-1}$~; en particulier, si $g \in \mathrm{fix}(\sigma)$, alors $\alpha_0(g)\in \mathrm{fix}(\sigma)$.

Soit maintenant $\sigma = (2 \ldots n-1) \in B_n$. Alors $\mathrm{fix}(\sigma)=\{0,[\sigma_{1,n}]\}$. Donc, puisque $\alpha_0([\sigma_{1,n}])\in \mathrm{fix}(\sigma)$, on a $\alpha_0([\sigma_{1,n}])=[\sigma_{1,n}]$. De même, pour tout $i \in \{1,\ldots,n-1\}$, $\alpha_0([\sigma_{i,n}])=[\sigma_{i,n}]$. Ainsi, $\alpha_0|_{F^{n-2}}=\mathrm{id}_{F^{n-2}}$. Puisque, par ailleurs, $\alpha_0$ est l'identité sur $B_n$, on voit que $\alpha_0|_{U_n}=\mathrm{id}_{U_n}$. De ce fait, étant donné que $\alpha_0|_{A_n}=\mathrm{id}_{A_n}$ et que $A_n$ et $U_n$ engendrent $\Out(W_n)$ par la proposition~\ref{système fini de générateurs}, on voit que $\alpha_0=\mathrm{id}$ et le résultat s'en déduit.

\subsection{Démonstration dans le cas $n=6$}

Dans le cas où $n=6$, la proposition~\ref{Sn fixe 0 étoile} s'appliquant encore, soit $\alpha_0$ un représentant de la classe d'automorphismes extérieurs de $\alpha$ tel que $\alpha_0(A_n)=A_n$. Supposons que la classe d'automorphisme extérieur de $\alpha_0|_{A_n}$ soit non triviale. Alors une description explicite d'un automorphisme engendrant l'unique classe d'automorphismes extérieurs de $\mathfrak{S}_6$ (cf. \cite{Miller58}) donne, en identifiant $A_n$ et $\mathfrak{S}_6$ par l'unique isomorphisme envoyant $\tau_i$ sur la permutation $(i \; i+1)$ pour $1 \leq i \leq 5$, que $$\alpha_0(B_n)=\left\langle [(1 \,2)(3 \,4)(5 \, 6)],[(1 \,6)(2 \,4)(3 \, 5)],[(1 \,4)(2 \,3)(5 \, 6)],[(1 \,6)(2 \,5)(3 \, 4)]\right\rangle.$$ Ainsi, $\alpha_0(B_n)$ agit transitivement sur les classes de conjugaison de $\{x_1,\ldots,x_n\}$. Alors, puisque $\alpha_0(B_n) \subseteq \alpha_0(U_n)$, par le quatrième point de la proposition~\ref{Wn fixe Z2 étoile}, $\alpha_0(B_n)$ fixe une $F$-étoile $\mathcal{X}$. Soit $X$ un représentant de $\mathcal{X}$. Par la proposition~\ref{prop Levitt}, le noyau du morphisme naturel $\alpha_0(B_n) \to \Aut_{gr}(X)$ est isomorphe à $F^{n-2} \cap \alpha_0(B_n)$.

Or $F^{n-2} \cap \alpha_0(B_n)$ est un $2$-sous-groupe distingué de $\alpha_0(B_n)$. Comme $\alpha_0(B_n)$ est isomorphe à $\mathfrak{S}_{n-1}$ et que $n=6$, nous avons $F^{n-2} \cap \alpha_0(B_n)=\{1\}$. Donc $\alpha_0(B_n)$ est isomorphe à $\Aut_{gr}(X)$ car $\Aut_{gr}(X)$ est isomorphe à $ \mathfrak{S}_{n-1}$. Soient maintenant $\overline{X}$ le graphe sous-jacent à $X$, $v_1,\ldots, v_{n-1}$ les feuilles de $\overline{X}$, et $v_n$ le centre de $\overline{X}$. Pour $j \in \{1,\ldots,n\}$, soit $\left\langle y_j \right\rangle$ l'image réciproque par le marquage du groupe associé à $v_j$. Le groupe $\Aut_{gr}(X)$, et donc $\alpha_0(B_n)$, s'identifie à l'ensemble des bijections de $\{v_1,\ldots,v_n\}$ fixant $v_n$. Or, par la proposition~\ref{système fini de générateurs}, il existe $\pi \in \mathrm{Bij}(\{x_1,\ldots,x_n\})$ telle que pour tout $i \in \{1,\ldots,n\}$, il existe $z_i \in W_n$ vérifiant~:
\begin{equation*}
y_i=z_ix_{\pi(i)}z_i^{-1}.
\end{equation*}
Ceci contredit le fait que $\alpha_0(B_n)$ s'identifie à l'ensemble des bijections de $\{v_1,\ldots,v_n\}$ fixant $v_n$ car le groupe $\alpha_0(B_n)$ agit transitivement sur l'ensemble des classes de conjugaison de $\{x_1,\ldots,x_n\}$. Donc la classe d'automorphisme extérieur de $\alpha_0|_{A_n}$ est triviale et on conclut comme dans~\ref{démonstration n=5}.

\subsection{Démonstration dans le cas $n=4$}

Dans le cas où $n=4$, la proposition~\ref{Sn fixe 0 étoile} et le quatrième point de la proposition~\ref{Wn fixe Z2 étoile} ne sont plus valables car alors tout sous-groupe de $\Out(W_n)$ isomorphe à $\mathfrak{S}_4$ est isomorphe au produit semi-direct $V \rtimes \mathfrak{S}_3$, où $V$ est le groupe de Klein. Nous avons cependant la proposition suivante.

\begin{prop}\label{Wn cas n=4}
Soient $n=4$ et $G$ un sous-groupe de $\Out(W_n)$ isomorphe au produit semi-direct $F^{n-2} \rtimes \mathfrak{S}_{n-1}$. Alors $G$ est soit le stabilisateur d'une unique $F$-étoile, soit le stabilisateur d'une unique $\{0\}$-étoile. Les deux cas sont mutuellement exclusifs.
\end{prop}

\dem Soient $\mathcal{X}$ un point de l'épine de $\PO(W_n)$ fixé par $G$ (qui existe par la proposition~\ref{point fixe outre espace}), et $X$ un représentant de $\mathcal{X}$. Soient $\overline{X}$ le graphe sous-jacent à $X$ et $L$ l'ensemble des feuilles de $\overline{X}$. La proposition~\ref{Wn cas n=4} se démontre de manière identique à la proposition~\ref{Wn fixe Z2 étoile}~$(3)$, à ceci près que l'on ne peut pas exclure le cas où $\overline{X}$ possède $n$ feuilles. Il faut alors distinguer le cas où $|L|=n-1$ et $|L|=n$. Si $\overline{X}$ possède $n$ feuilles, le lemme~\ref{Sn n feuilles} donne que $\mathcal{X}$ est une $\{0\}$-étoile. Si $\overline{X}$ possède $n-1$ feuilles, alors la proposition~\ref{Wn fixe Z2 étoile}~$(3)$ donne que $\mathcal{X}$ est une $F$-étoile.

\bigskip

Montrons maintenant que $G$ ne peut fixer à la fois une $\{0\}$-étoile et une $F$-étoile. Par la proposition~\ref{Wn fixe Z2 étoile}~$(1)$, $G$ est le stabilisateur de tout point fixé par $G$. 

Supposons que $G$ soit le stabilisateur d'une $\{0\}$-étoile $\mathcal{X}$. Soient $X$ un représentant de $\mathcal{X}$, et $\overline{X}$ le graphe sous-jacent à $X$. Soient $v_1,\ldots,v_n$ les sommets de $\overline{X}$ dont les groupes associés sont non triviaux et, pour tout $i \in \{1,\ldots,n\}$, soit $y_i$ l'image réciproque par le marquage du générateur du groupe associé à $v_i$. Alors le groupe $G$ est le groupe engendré par les permutations de $\{y_1,\ldots,y_n\}$. 

Soit $\mathcal{Y}$ une $F$-étoile dans l'épine de $\PO(W_n)$ fixée par $G$. Par le corollaire~\ref{points fixes connexes}, $\mathrm{Fix}_{K_n}(G)$ est connexe. Il existe donc un chemin continu dans $\mathrm{Fix}_{K_n}(G)$ de $\mathcal{X}$ vers $\mathcal{Y}$. Les sommets par lesquels passe ce chemin sont uniquement des $\{0\}$-étoiles et des $F$-étoiles au vu des points stabilisés par $G$. Or le groupe engendré par les permutations de $\{y_1,\ldots,y_n\}$ ne fixe aucune $F$-étoile adjacente à $X$. En effet, le groupe $G$ contiendrait un élément permutant le centre de la $F$-étoile avec une feuille, ce qui n'est pas possible. Donc $G$ ne fixe aucune $F$-étoile.

\bigskip

Enfin, l'unicité du point fixe provient du fait que l'ensemble des $\{0\}$-étoiles et l'ensem\-ble des $F$-étoiles sont discrets dans l'épine de $\PO(W_n)$ alors que l'ensemble des points fixes de $G$ est connexe par le corollaire~\ref{points fixes connexes}. \hfill\qedsymbol

\bigskip

Nous pouvons maintenant montrer le théorème~\ref{Théorème out coxeter} dans le cas $n=4$.

Soit $\alpha \in \Aut(\Out(W_n))$. Soit $\mathcal{X}_1$ la $\{0\}$-étoile fixée par le sous-groupe fini $A_n \simeq \mathfrak{S}_4$ de $\Out(W_n)$. Par la proposition~\ref{Wn cas n=4}, $\alpha(A_n)$ fixe soit une $\{0\}$-étoile, soit une $F$-étoile. 

Si $\alpha(A_n)$ fixe une $\{0\}$-étoile, alors la même démonstration que pour le cas où \mbox{$n \neq 6$} dans la partie~\ref{démonstration n=5} montre que quitte à changer $\alpha$ dans sa classe d'automorphisme extérieurs, nous avons $\alpha|_{A_n}=\mathrm{id}_{A_n}$. Par la proposition~\ref{Wn cas n=4}, le groupe $U_n \simeq F^2 \rtimes \mathfrak{S}_3$ fixe soit une $\{0\}$-étoile, soit une $F$-étoile. Étant donné que $B_n \subseteq U_n$ fixe une unique $\{0\}$-étoile $\rho$ et une unique $F$-étoile $\rho'$ et que $\alpha|_{B_n}=\mathrm{id}_{B_n}$, on voit que $\alpha(U_n)$ est soit le stabilisateur de $\rho$, soit le stabilisateur de $\rho'$. Cependant, puisque le stabilisateur de $\rho$ est $A_n$ et que $\alpha|_{A_n}=\mathrm{id}_{A_n}$, on voit que $\alpha(U_n)$ est le stabilisateur de $\rho'$. Donc $\alpha(U_n)=U_n$. Le reste de la démonstration est alors identique à celle du cas où $n \neq 6$ dans la partie~\ref{démonstration n=5}.

Supposons que $\alpha(A_n)$ fixe une unique $F$-étoile. Construisons à présent un représentant de la classe d'automorphismes extérieurs de $\alpha$. Puisque $\Out(W_n)$ agit transitivement sur les $F$-étoiles, quitte à changer $\alpha$ dans sa classe d'automorphismes extérieurs, on peut supposer que $\alpha(A_n)=U_n$. Soit $V$ le groupe de Klein contenu dans $A_n$. Alors $\alpha(V)$ est l'unique $2$-sous-groupe distingué non trivial de $U_n$. Donc $$\alpha(V)=\left\langle [\sigma_{1,4}],[\sigma_{2,4}], [\sigma_{3,4}] \right\rangle.$$
Ainsi, puisque $B_n \cap V= \{\mathrm{id}\}$, on voit que $\alpha(B_n)\cap \alpha(V)=\{\mathrm{id}\}$. Par ailleurs, $A_n=B_nV$, donc $U_n=\alpha(B_n)\alpha(V)$. De ce fait, $\alpha(B_n)$ est un sous-groupe de $U_n$ d'ordre $6$. Or, il existe une unique classe de conjugaison de sous-groupes d'ordre $6$ dans $U_n$. Donc, quitte à changer $\alpha$ dans sa classe d'automorphismes extérieurs, on peut supposer que $\alpha(B_n)=B_n$. De même, puisque $B_n$ est isomorphe à $\mathfrak{S}_3$, quitte à changer $\alpha$ dans sa classe d'automorphisme extérieur, on peut supposer que $\alpha|_{B_n}=\mathrm{id}_{B_n}$.

Déterminons à présent l'image de $[\tau_3]$ et $[\sigma_{3,4}]$ par $\alpha$. Puisque $[\tau_1][\tau_3] \in V$, on voit que $\alpha([\tau_1][\tau_3]) \in \{[\sigma_{1,4}],[\sigma_{2,4}], [\sigma_{3,4}]\}$. Or, $[\tau_1]$ commute avec $[\tau_1][\tau_3]$, donc $\alpha([\tau_1][\tau_3])$ doit également commuter avec $[\tau_1]$. De ce fait, $\alpha([\tau_1][\tau_3])=[\sigma_{3,4}]$ et $\alpha([\tau_3])=[\tau_1][\sigma_{3,4}]$.

Déterminons l'image de $[\sigma_{3,4}]$ par $\alpha$. Puisque $\alpha(B_n)=B_n$, le groupe $\alpha(U_n)$ est le stabilisateur d'un point fixe de $B_n$. Par la proposition~\ref{unicité Z2 étoile}, $B_n$ fixe uniquement deux sommets de l'épine de $\PO(W_n)$~: la $\{0\}$-étoile stabilisée par $A_n$ et la $F$-étoile stabilisée par $U_n$. Comme $\alpha(A_n)=U_n$, on a nécessairement $\alpha(U_n)=A_n$. Donc $\alpha([\sigma_{3,4}]) \in V$. Puisque $[\sigma_{3,4}]$ commute avec $[\tau_1]$, on obtient que $\alpha([\sigma_{3,4}])=[\tau_1][\tau_3]$.

\bigskip

Donc $\alpha$ se restreint en l'identité sur $B_n$, envoie $[\tau_3]$ sur $[\tau_1][\sigma_{3,4}]$ et $[\sigma_{3,4}]$ sur $[\tau_1][\tau_3]$. Comme $B_n$, $[\tau_3]$ et $[\sigma_{3,4}]$ engendrent $\Out(W_4)$, ceci montre qu'un tel automorphisme $\alpha$, s'il existe, est unique modulo automorphisme intérieur.

Réciproquement, montrons que l'application $\alpha$ de $B_n \cup \{[\tau_3],[\sigma_{3,4}]\}$ dans $\Out(W_4)$ définie par $\alpha|_{B_n}=\mathrm{id}_{B_n}$, $\alpha([\tau_3])=[\tau_1][\sigma_{3,4}]$ et $\alpha([\sigma_{3,4}])=[\tau_1][\tau_3]$ s'étend de manière unique en un morphisme de groupes de $\Out(W_4)$. Comme $[\tau_1]$ commute avec $[\tau_3]$ et $[\sigma_{3,4}]$, ceci montre que $\alpha$ est involutif, donc un automorphisme de $\Out(W_4)$. Sa classe dans $\Out(\Out(W_4))$ est non triviale (car son action sur l'épine de $\PO(W_4)$ est non triviale), ce qui montre le théorème~\ref{Théorème out coxeter} lorsque $n=4$.

Pour simplifier les notations, nous notons $[i \; j]$ la classe d'automorphismes extérieurs de la transposition permutant $x_i$ et $x_j$. Notons $$S=\{[i\;j] \;|\; 1 \leq i,j \leq 4\} \cup \{[\sigma_{i,j}]\;|\;1\leq i \neq j \leq 4\},$$ qui est une partie génératrice de $\Out(W_4)$ par la proposition~\ref{système fini de générateurs}. Un petit calcul élémentaire montre que, si $i=1,2$, alors 
$$
[i\;4]=[i\;3][3\;4][i\;3], [\sigma_{i,4}]=[i\;3][ \sigma_{3,4}][i\;3],$$ 
$$\alpha([i\;3])\alpha([3\;4])\alpha([i\;3])=[j\;k][ \sigma_{i,4}]\text{ et } \alpha([i\;3])\alpha([ \sigma_{3,4}])\alpha([i\;3])=[j\;k][i\;4],$$

où $\{j,k\}=\{1,2,3\}-\{i\}$. Considérons l'application $\widetilde{\alpha}$ de $S$ dans $\Out(W_4)$ étendant $\alpha$ sur $S\cap (B_n \cup \{[3\;4],[\sigma_{3,4}]\})$ et telle que, si $i=1,2$, $$\widetilde{\alpha}([i\;4])=[j\;k][\sigma_{i,4}] \text{ et } \widetilde{\alpha}([\sigma_{i,4}])=[j\;k][i\;4],$$ où $\{j,k\}=\{1,2,3\}-\{i\}$.
Des calculs élémentaires pour lesquels nous renvoyons à l'appendice~\ref{appendice} montrent que cette application préserve, quand $n=4$, la présentation de $\Out(W_n)$ donnée par~\cite[Theorem 2.20]{Gilbert1987}, ce qui conclut.

\subsection{Démonstration de la rigidité de $\Aut(W_n)$}

Nous démontrons à présent le théorème~\ref{Théorème aut coxeter}. Soient $n \geq 4$ et $\alpha \in \Aut(\Aut(W_n))$. Soient $\widetilde{A}_n=\left\langle \tau_1,\ldots, \tau_{n-1} \right\rangle$, $\widetilde{B}_n=\left\langle \tau_1,\ldots, \tau_{n-2} \right\rangle$ et $\widetilde{U}_n=\left\langle \tau_1,\ldots, \tau_{n-2},\sigma_{1,n} \right\rangle$. En utilisant les remarques~\ref{Sn fixe 0étoile aut},~\ref{Wn fixe Z2 étoile Aut}~$(4)$ et~\ref{Unicité Z2 étoile Aut}~$(2)$, et en effectuant une démonstration identique à celle du théorème~\ref{Théorème out coxeter} dans les cas où $n \geq 5$, on voit que, quitte à changer $\alpha$ dans sa classe d'automorphismes extérieurs, $\alpha|_{\widetilde{A}_n}=\mathrm{id}_{\widetilde{A}_n}$ et que $\alpha(\widetilde{U}_n)=\widetilde{U}_n$. 

Or $\widetilde{U}_n$ est isomorphe à $F^{n-1} \rtimes \widetilde{B}_n$. Soit $\sigma \in \widetilde{B}_n$. On note $\mathrm{fix}(\sigma)$ l'ensemble des points fixes de $\sigma$ agissant par conjugaison dans $F^{n-1}$. On voit que pour tout $\sigma \in \{0\} \rtimes \widetilde{B}_n$ et pour tout $g \in F^{n-1} \rtimes \{1\}$, $\alpha(\sigma g \sigma^{-1})=\sigma\alpha(g)\sigma^{-1}$~; en particulier, si $g \in \mathrm{fix}(\sigma)$, alors $\alpha(g)\in \mathrm{fix}(\sigma)$.

\bigskip

Soit maintenant $\sigma = (2 \ldots n-1) \in B_n$. Alors $\mathrm{fix}(\sigma)=\{0,\sigma_{1,n}, \prod\limits_{i\neq 1,n} \sigma_{i,n}, \prod\limits_{i=1}^{n-1} \sigma_{i,n}\}$. Donc $\alpha(\sigma_{1,n}) \in \{\sigma_{1,n}, \prod\limits_{i\neq 1,n} \sigma_{i,n}, \prod\limits_{i=1}^{n-1} \sigma_{i,n}\}$. Comme $\prod\limits_{i=1}^{n-1} \sigma_{i,n}$ est l'unique élément non trivial dans le centre de $\widetilde{U}_n$, on voit que $\alpha(\sigma_{1,n}) \neq \prod\limits_{i=1}^{n-1} \sigma_{i,n}$.

Supposons par l'absurde que $\alpha(\sigma_{1,n})=\prod\limits_{i\neq 1,n} \sigma_{i,n}$. Pour $j \in \{1,\ldots,n-1\}$, notons $(1 \;j)$ la transposition de $\widetilde{B}_n$ permutant $x_1$ et $x_j$. Alors, on voit que, pour tout $j \in \{1,\ldots,n-1\}$, $\alpha(\sigma_{j,n})= \alpha((1 \;j)\sigma_{1,n}(1\;j))=\prod\limits_{i\neq j,n} \sigma_{i,n}$.

Un calcul immédiat montre alors que, pour tout $j \neq k,n$, et $k<n$,

$$\alpha(\sigma_{k,j})= \alpha((j \;n)\sigma_{k,n}(j\;n))=\prod\limits_{i\neq j,k} \sigma_{i,j}.$$

Or $\sigma_{1,2}\sigma_{3,4}=\sigma_{3,4}\sigma_{1,2}$, alors que

$$\begin{array}{cc}
\alpha(\sigma_{1,2})\alpha(\sigma_{3,4})(x_1)=
\prod\limits_{i\neq 1,2} \sigma_{i,2}\prod\limits_{i\neq 3,4} \sigma_{i,4}(x_1)=x_2x_4x_2x_1x_2x_4x_2 & \text{ et que } \\
\alpha(\sigma_{3,4})\alpha(\sigma_{1,2})(x_1)=
\prod\limits_{i\neq 3,4} \sigma_{i,4}\prod\limits_{i\neq 1,2} \sigma_{i,2}(x_1)=x_4x_1x_4.
\end{array}
$$

Donc $\alpha(\sigma_{1,2})\alpha(\sigma_{3,4}) \neq \alpha(\sigma_{3,4})\alpha(\sigma_{1,2})$. Ceci contredit le fait que $\alpha$ est un morphisme de groupes. Ainsi, $\alpha(\sigma_{1,n})=\sigma_{1,n}$. Par la proposition~\ref{système fini de générateurs}, nous avons $\alpha=\mathrm{id}$. Ceci conclut la démonstration du théorème~\ref{Théorème aut coxeter}.

\appendix
\section{Présentation du groupe $\Out(W_4)$}\label{appendice}

Soit $n=4$. Pour simplifier les notations, nous notons $[i \; j]$ la classe d'automorphismes extérieurs de la transposition permutant $x_i$ et $x_j$. Nous rappelons que l'application ensembliste ${\alpha \colon B_n \cup \{[\tau_3],[\sigma_{3,4}]\} \to \Out(W_n)}$ est définie par~:
$$\begin{array}{c}
\alpha|_{B}=\mathrm{id}_{B_n}~; \\
\alpha([3\;4])=[1\;2][\sigma_{3,4}]~; \\
\alpha([\sigma_{3,4}])=[1\;2][3\;4].
\end{array}$$

Nous montrons dans cette appendice que l'application $\alpha$ s'étend de manière unique en un morphisme de groupes de $\Out(W_4)$ dans lui-même. Nous montrons pour cela qu'il préserve l'ensemble des relations d'une présentation de $\Out(W_4)$. La présentation suivante est due à Gilbert.

\begin{prop}\cite[Theorem 2.20]{Gilbert1987}\label{présentation AutWn}
Soit $n$ un entier plus grand que $2$. Une présentation de $\Out(W_n)$ est donnée par~: 
\begin{enumerate}
\item la partie génératrice $S$ constituée de l'ensemble des permutations $[i\;j]$ pour les entiers distincts $i,j \in \{1,\ldots,n\}$ ainsi que l'ensemble des éléments $[\sigma_{i,j}]$ pour les entiers distincts $i,j \in \{1,\ldots,n\}$~;
\item les relations suivantes~:
\begin{enumerate}
\item pour tout $i \in \{1,\ldots,n\}$, $\prod_{j \neq i}=[\sigma_{j,i}]=1$~;
\item pour tous les $i,j,k,\ell \in \{1,\ldots,n\}$ avec $i \neq j$ et $k \neq \ell$, si on pose $\tau=(i\;j)$, alors $[i \;j][k\;\ell]=[\tau(k)\;\tau(\ell)][i \;j]$~;
\item pour tout $j \in \{1,\ldots,n\}$, pour tous les $i,k \in \{1,\ldots,n\}-\{j\}$ distincts, $[\sigma_{i,j}][\sigma_{k,j}]=[\sigma_{k,j}][\sigma_{i,j}]$~;
\item pour tous les $i,j  \in \{1,\ldots,n\}$ distincts, $[\sigma_{i,j}][\sigma_{i,j}]=1$~;
\item pour tous les $i,j,k,\ell \in \{1,\ldots,n\}$ deux à deux distincts, $[\sigma_{i,j}][\sigma_{k,\ell}]=[\sigma_{k,\ell}][\sigma_{i,j}]$~;
\item pour tous les $i,j,k,\ell \in \{1,\ldots,n\}$, tels que $k \neq \ell$, si $\tau=(i \; j)$, $[i\;j][\sigma_{k,\ell}]=[\sigma_{\tau(k),\tau(\ell)}][i\;j]$~;
\item pour tous les $i,j,k \in \{1,\ldots,n\}$ deux à deux distincts, $[\sigma_{j,i}][\sigma_{i,k}][\sigma_{j,k}]=[\sigma_{j,k}][\sigma_{i,k}][\sigma_{j,i}]$.
\end{enumerate}
\end{enumerate}
\end{prop}

Nous remarquons que, dans le cas où $n=4$, la relation $(g)$ se déduit des relations $(a)$, $(d)$ et $(e)$.

\begin{prop}
L'application $\alpha$ se prolonge de manière unique en un morphisme de groupes de $\Out(W_4)$ dans lui-même.
\end{prop}

\dem Nous définissons tout d'abord une application prolongeant $\alpha$ sur la partie génératrice $S$ de $\Out(W_4)$ définie dans la proposition~\ref{présentation AutWn}. Un petit calcul élémentaire montre que, si $i=1,2$, $[i\;4]=[i\;3][3\;4][i\;3]$, $[\sigma_{i,4}]=[i\;3][\sigma_{3,4}][i\;3]$, $\alpha([i\;3])\alpha([3\;4])\alpha([i\;3])=[j\;k][\sigma_{i,4}]$, $\alpha([i\;3])\alpha([\sigma_{3,4}])\alpha([i\;3])=[j\;k][i\;4]$, où $\{j,k\}=\{1,2,3\}-\{i\}$.

Nous considérons à présent l'application $\widetilde{\alpha}$ de $S$ dans $\Out(W_4)$ étendant $\alpha$ sur $S\cap (B_n \cup \{[3\;4],[\sigma_{3,4}]\})$ et telle que, si $i=1,2$, $\widetilde{\alpha}([i\;4])=[j\;k][\sigma_{i,4}]$, $\widetilde{\alpha}([\sigma_{i,4}])=[j\;k][i\;4]$, où $\{j,k\}=\{1,2,3\}-\{i\}$ et, si $i$ et $j$ sont distincts et si $j \neq 4$, $\widetilde{\alpha}([\sigma_{i,j}])=[\sigma_{j,4}][i\;j][k\;\ell][\sigma_{j,4}]$, où $\{k,\ell\}=\{1,2,3,4\}-\{i,j\}$. 

Vérifions maintenant que $\widetilde{\alpha}$ préserve la présentation de $\Aut(W_n)$. Ceci montrera que $\widetilde{\alpha}$ se prolonge en un morphisme de groupes de $\Out(W_4)$ dans lui-même. De plus, étant donné que $B \cup \{[3\;4],[\sigma_{3,4}]\}$ est une partie génératrice de $\Out(W_4)$ (cf. \cite[Theorem B]{muhlherr1997}), au vu de la définition de $\widetilde{\alpha}$, ceci conclura la démonstration de la proposition. Nous écrivons pour chaque cas la relation vérifiée en préalable à la démonstration.

\bigskip\noindent{$(1)$ } \textit{Pour tout $i$, pour tous les $j,k,\ell\in \{1,2,3,4\}-\{i\}$ deux à deux distincts, $[\sigma_{j,i}][\sigma_{k,i}][\sigma_{\ell,i}]=1$}.

$$\begin{array}{c}
\widetilde{\alpha}([\sigma_{1,4}])\widetilde{\alpha}([\sigma_{2,4}])\widetilde{\alpha}([\sigma_{3,4}])=[2\;3][1\;4][1\;3][2\;4][1\;2][3\;4]=1.
\end{array}$$
Si $j \neq 4$, et si $i,k,\ell \in \{1,2,3,4\}-\{j\}$ sont deux à deux distincts,
$$\begin{array}{ccl}
\widetilde{\alpha}([\sigma_{i,j}])\widetilde{\alpha}([\sigma_{k,j}])\widetilde{\alpha}([\sigma_{\ell,j}])&=&[\sigma_{j,4}][i\;j][k\;\ell][\sigma_{j,4}][\sigma_{j,4}][j\;k][i\;\ell][\sigma_{j,4}][\sigma_{j,4}][j\;\ell][i\;k][\sigma_{j,4}] \\
{} &=& [\sigma_{j,4}][i\;j][k\;\ell][j\;k][i\;\ell][j\;\ell][i\;k][\sigma_{j,4}]=[\sigma_{j,4}][\sigma_{j,4}]=1.
\end{array}$$

\noindent{$(2)$ } \textit{Pour tous les $i,j,k,\ell$ vérifiant $i \neq j$ et $k \neq \ell$, si on pose $\sigma=(i\;j)$, alors $[i \;j][k\;\ell]=[\sigma(k)\;\sigma(\ell)][i \;j]$.}

\bigskip

Puisque $\widetilde{\alpha}$ est l'identité sur $B_n$, cette relation est vérifiée si $i,j,k,\ell \in \{1,2,3\}$. Vérifions les autres cas. Soient $i,j \in \{1,2,3\}$ distincts et $k \in \{1,2,3\}-\{i,j\}$.
$$\begin{array}{cl}
\widetilde{\alpha}([i\;4])\widetilde{\alpha}([j \;4])&=[j\;k][\sigma_{i,4}][i\;k][\sigma_{j,4}]=[j\;i][j\;k][\sigma_{k,4}][\sigma_{j,4}] \\
{} &=[j\;i][j\;k][\sigma_{i,4}]=\widetilde{\alpha}([j\;i])\widetilde{\alpha}([i\;4]).
\end{array}
$$
Maintenant, si $i,j,k \in \{1,2,3\}$ sont deux à deux distincts,

$$\widetilde{\alpha}([i\;j])\widetilde{\alpha}([k \;4])=[i\;j][i\;j][\sigma_{k,4}]=[\sigma_{k,4}]=\widetilde{\alpha}([k\;4])\widetilde{\alpha}([i\;j]).$$

\bigskip\noindent{$(3)$ } \textit{Pour tout $j$, pour tous les $i,k \in \{1,2,3,4\}-\{j\}$ distincts, nous avons $[\sigma_{i,j}][\sigma_{k,j}]=[\sigma_{k,j}][\sigma_{i,j}]$.} 

\bigskip

On note $\ell$ l'élément distinct de $i,j$ et $k$.

Supposons que $j \neq 4$. Alors
$$\begin{array}{ccl}
\widetilde{\alpha}([\sigma_{i,j}])\widetilde{\alpha}([\sigma_{k,j}])&=& [\sigma_{j,4}][j\;i][k\;\ell][j\;k][i\;\ell][\sigma_{j,4}] \\
{} & = & [\sigma_{j,4}][j\;\ell][i\;k][\sigma_{j,4}]=\widetilde{\alpha}([\sigma_{\ell,j}]) \\
{} &=& \widetilde{\alpha}([\sigma_{k,j}])\widetilde{\alpha}([\sigma_{i,j}]).
\end{array}$$

Dans le cas où $j=4$,
$$\begin{array}{ccl}
\widetilde{\alpha}([\sigma_{i,4}])\widetilde{\alpha}([\sigma_{k,4}])&=& [k\;\ell][i\;4][i\;\ell][k\;4]=[i\;\ell][k\;4][\ell\;k][i\;4]=\widetilde{\alpha}([\sigma_{k,4}])\widetilde{\alpha}([\sigma_{i,4}]).
\end{array}$$

\bigskip\noindent{$(4)$ } \textit{Pour tout $i \neq j$, nous avons $[\sigma_{i,j}][\sigma_{i,j}]=1$. }

\bigskip

On note $\ell$ et $k$ les deux éléments distincts de $i$ et $j$.

Supposons que $j \neq 4$. Alors
$$\begin{array}{ccl}
\widetilde{\alpha}([\sigma_{i,j}])\widetilde{\alpha}([\sigma_{i,j}])&=& [\sigma_{j,4}][j\;i][k\;\ell][j\;i][k\;\ell][\sigma_{j,4}]=1.
\end{array}$$

Si $j=4$, alors
$$\begin{array}{ccl}
\widetilde{\alpha}([\sigma_{i,4}])\widetilde{\alpha}([\sigma_{i,4}])&=& [k\;\ell][i\;4][k\;\ell][i\;4]=1.
\end{array}$$

\bigskip\noindent{$(5)$ } \textit{Si $i,j,k,\ell$ sont deux à deux distincts, alors $[\sigma_{i,j}][\sigma_{k,\ell}]=[\sigma_{k,\ell}][\sigma_{i,j}]$.}

\bigskip

Nous traitons tout d'abord le cas où $j=4$ ou $\ell=4$. Par symétrie, nous pouvons supposer que $j=4$.

$$\begin{array}{ccl}
\widetilde{\alpha}([\sigma_{i,4}])\widetilde{\alpha}([\sigma_{k,\ell}])&=&[k\;\ell][i\;4][\sigma_{\ell,4}][k\;\ell][i\;4][\sigma_{\ell,4}] \\
{} &=& [k\;\ell][\sigma_{\ell,i}][k\;\ell][i\;4][\sigma_{\ell,i}][i\;4] \\
{} &=& [\sigma_{k,i}][k\;\ell][i\;4][\sigma_{k,i}][k\;\ell][i\;4] \\
{} &=& [\sigma_{\ell,4}][\sigma_{k,i}][\sigma_{\ell,4}][k\;\ell][i\;4][\sigma_{\ell,4}][\sigma_{k,i}][\sigma_{\ell,4}][k\;\ell][i\;4] \\
{} &=& [\sigma_{\ell,4}][k\;\ell][\sigma_{\ell,i}][\sigma_{k,4}][\sigma_{\ell,i}][\sigma_{k,4}][i\;4][\sigma_{\ell,4}][k\;\ell][i\;4] \\
{} &=& [\sigma_{\ell,4}][k\;\ell][i\;4][\sigma_{\ell,4}][k\;\ell][i\;4]=\widetilde{\alpha}([\sigma_{k,\ell}])\widetilde{\alpha}([\sigma_{i,4}]).
\end{array}$$
Nous effectuons maintenant le cas où $i=4$ ou $k=4$. Par symétrie, nous pouvons supposer $i=4$.
$$\begin{array}{ccl}
\widetilde{\alpha}([\sigma_{4,j}])\widetilde{\alpha}([\sigma_{k,\ell}])&=& [\sigma_{j,4}][j\;4][k\;\ell][\sigma_{j,4}][\sigma_{\ell,4}][k\;\ell][j\;4][\sigma_{\ell,4}] \\
{} &=& [\sigma_{j,4}][j\;4][k\;\ell][\sigma_{k,4}][k\;\ell][j\;4][\sigma_{\ell,4}] \\
{} &=& [\sigma_{j,4}][\sigma_{\ell,j}][\sigma_{\ell,4}]~; \\
{} & {} & {} \\
\widetilde{\alpha}([\sigma_{k,\ell}])\widetilde{\alpha}([\sigma_{4,j}])&=& [\sigma_{\ell,4}][j\;4][k\;\ell][\sigma_{\ell,4}][\sigma_{j,4}][k\;\ell][j\;4][\sigma_{j,4}] \\
{} &=& [\sigma_{\ell,4}][j\;4][k\;\ell][\sigma_{k,4}][k\;\ell][j\;4][\sigma_{j,4}] \\
{} &=& [\sigma_{\ell,4}][\sigma_{\ell,j}][\sigma_{j,4}] \\
{} &=& [\sigma_{j,4}][\sigma_{\ell,4}][\sigma_{j,4}][\sigma_{\ell,j}][\sigma_{\ell,4}][\sigma_{j,4}][\sigma_{\ell,4}] \\
{} &=& [\sigma_{j,4}][\sigma_{k,4}][\sigma_{\ell,j}][\sigma_{k,4}][\sigma_{\ell,4}] \\
{} &=& [\sigma_{j,4}][\sigma_{\ell,j}][\sigma_{\ell,4}]=\widetilde{\alpha}([\sigma_{4,j}])\widetilde{\alpha}([\sigma_{k,\ell}]).
\end{array}$$

\bigskip\noindent{$(6)$ } \textit{Pour tous les $i,j,k,\ell$ tels que $k \neq \ell$, si $\tau=(i \; j)$, alors $[i\;j][\sigma_{k,\ell}]=[\sigma_{\tau(k),\tau(\ell)}][i\;j]$.} 

\bigskip

On note $a$ et $b$ les éléments vérifiant $\{a,b\}=\{1,2,3,4\}-\{k,\ell\}$.

Nous supposons tout d'abord que $i,j \in \{1,2,3\}$. Supposons également que $\ell \neq 4$. Si $\{i,j\} \cap \{k,\ell\}=\varnothing$, alors $\{a,b\}=\{i,j\}$ et $\tau(k)=k$ et $\tau(\ell)=\ell$. Donc

$$\begin{array}{ccl}
\widetilde{\alpha}([i\;j])\widetilde{\alpha}([\sigma_{k,\ell}])&=&[i\;j][\sigma_{\ell,4}][k\;\ell][i\;j][\sigma_{\ell,4}] =[\sigma_{\ell,4}][k\;\ell][i\;j][\sigma_{\ell,4}][i\;j]=\widetilde{\alpha}([\sigma_{k,\ell}])\widetilde{\alpha}([i\;j]).
\end{array}$$

Si $\{i,j\} \cap \{k,\ell\}=\{k\}=\{i\}$, alors $\{a,b\}=\{j,4\}$ et $\tau(k)=j$ et $\tau(\ell)=\ell$. Donc
$$\begin{array}{ccl}
\widetilde{\alpha}([i\;j])\widetilde{\alpha}([\sigma_{i,\ell}])&=&[i\;j][\sigma_{\ell,4}][i\;\ell][j\;4][\sigma_{\ell,4}] =[\sigma_{\ell,4}][j\;\ell][i\;4][\sigma_{\ell,4}][i\;j]=\widetilde{\alpha}([\sigma_{j,\ell}])\widetilde{\alpha}([i\;j]).
\end{array}$$

Si $\{i,j\} \cap \{k,\ell\}=\{\ell\}=\{i\}$, alors $\{a,b\}=\{j,b\}$ avec $b \notin \{i,k\}$ et $\tau(k)=k$ et $\tau(\ell)=j$. Donc
$$\begin{array}{ccl}
\widetilde{\alpha}([i\;j])\widetilde{\alpha}([\sigma_{k,i}])&=&[i\;j][\sigma_{i,4}][k\;i][j\;b][\sigma_{i,4}] =[\sigma_{j,4}][k\;j][i\;b][\sigma_{j,4}][i\;j]=\widetilde{\alpha}([\sigma_{k,j}])\widetilde{\alpha}([i\;j]).
\end{array}$$

Si $\{i,j\} \cap \{k,\ell\}=\{i,j\}$, alors $\{i,j\} \cap \{a,b\}=\varnothing$. De plus, puisque $i$ et $j$ jouent un rôle symétrique, nous pouvons supposer que $\tau(k)=j$ et $\tau(\ell)=i$. Donc
$$\begin{array}{ccl}
\widetilde{\alpha}([i\;j])\widetilde{\alpha}([\sigma_{i,j}])&=&[i\;j][\sigma_{j,4}][i\;j][a\;b][\sigma_{j,4}] =[\sigma_{i,4}][i\;j][a\;b][\sigma_{i,4}][i\;j]=\widetilde{\alpha}([\sigma_{j,i}])\widetilde{\alpha}([i\;j]).
\end{array}$$

Supposons maintenant que $\ell=4$. 

Si $\{i,j\} \cap \{k,\ell\}=\varnothing$, alors $\{a,b\}=\{i,j\}$ et $\tau(k)=k$ et $\tau(\ell)=\ell$. Donc
$$\begin{array}{ccl}
\widetilde{\alpha}([i\;j])\widetilde{\alpha}([\sigma_{k,4}])&=&[i\;j][i\;j][k\;4]=[i\;j][k\;4][i\;j]=\widetilde{\alpha}([\sigma_{k,4}])\widetilde{\alpha}([i\;j]).
\end{array}$$

Si $\{i,j\} \cap \{k,\ell\}=\{k\}=\{i\}$, alors $\{a,b\}=\{j,b\}$, avec $b\neq i$ et $b \neq 4$ et $\tau(k)=j$ et $\tau(\ell)=\ell$. Donc
$$\begin{array}{ccl}
\widetilde{\alpha}([i\;j])\widetilde{\alpha}([\sigma_{i,4}])&=&[i\;j][j\;b][i\;4] =[i\;b][j\;4][i\;j]=\widetilde{\alpha}([\sigma_{j,4}])\widetilde{\alpha}([i\;j]).
\end{array}$$

Supposons maintenant que $j=4$. Supposons également que $\ell \neq 4$. Puisque $j=4$, le cas où $\ell=4$ est symétrique au cas où $k=4$.

Si $\{i,4\} \cap \{k,\ell\}=\varnothing$, alors $\{a,b\}=\{i,4\}$ et $\tau(k)=k$ et $\tau(\ell)=\ell$. Donc

$$\begin{array}{ccl}
\widetilde{\alpha}([i\;4])\widetilde{\alpha}([\sigma_{k,\ell}])&=&[k\;\ell][\sigma_{i,4}][\sigma_{\ell,4}][k\;\ell][i\;4][\sigma_{\ell,4}] \\
{} &=& [\sigma_{i,4}][\sigma_{k,4}][k\;\ell][k\;\ell][i\;4][\sigma_{\ell,4}] \\
{} &=& [\sigma_{\ell,4}][k\;\ell][i\;4][\sigma_{k,4}][k\;\ell] \\
{} &=& [\sigma_{\ell,4}][k\;\ell][i\;4][\sigma_{i,4}][\sigma_{k,4}][\sigma_{i,4}][k\;\ell]= \\
{} &=& [\sigma_{\ell,4}][k\;\ell][i\;4][\sigma_{\ell,4}][\sigma_{i,4}][k\;\ell]= \widetilde{\alpha}([\sigma_{k,\ell}])\widetilde{\alpha}([i\;4]).
\end{array}$$

Si $\{i,4\} \cap \{k,\ell\}=\{k\}=\{i\}$, alors $\{a,b\}=\{a,4\}$, avec $a \neq i$ et $a \neq \ell$ et $\tau(k)=4$ et $\tau(\ell)=\ell$. Donc

$$\begin{array}{ccl}
\widetilde{\alpha}([i\;4])\widetilde{\alpha}([\sigma_{i,\ell}])&=&[a\;\ell][\sigma_{i,4}][\sigma_{\ell,4}][i\;\ell][a\;4][\sigma_{\ell,4}] \\
{} &=& [\sigma_{i,4}][\sigma_{a,4}][a\;\ell][i\;\ell][a\;4][\sigma_{\ell,4}] \\
{} &=& [\sigma_{\ell,4}][a\;\ell][i\;\ell][a\;4][\sigma_{\ell,4}] \\
{} &=& [\sigma_{\ell,4}][i\;a][\ell\;4][\sigma_{a,4}][a\;\ell] \\
{} &=& [\sigma_{\ell,4}][i\;a][\ell\;4][\sigma_{i,4}][\sigma_{a,4}][\sigma_{i,4}][a\;\ell] \\
{} &=& [\sigma_{\ell,4}][i\;a][\ell\;4][\sigma_{\ell,4}][\sigma_{i,4}][a\;\ell]= \widetilde{\alpha}([\sigma_{4,\ell}])\widetilde{\alpha}([i\;4]).
\end{array}$$

Si $\{i,4\} \cap \{k,\ell\}=\{k\}=\{4\}$, alors $\{a,b\}=\{a,i\}$, avec $a \neq 4$ et $a \neq \ell$ et $\tau(k)=i$ et $\tau(\ell)=\ell$. Donc

$$\begin{array}{ccl}
\widetilde{\alpha}([i\;4])\widetilde{\alpha}([\sigma_{4,\ell}])&=&[a\;\ell][\sigma_{i,4}][\sigma_{\ell,4}][4\;\ell][a\;i][\sigma_{\ell,4}] \\
{} &=& [\sigma_{i,4}][\sigma_{a,4}][a\;\ell][4\;\ell][a\;i][\sigma_{\ell,4}] \\
{} &=& [\sigma_{\ell,4}][a\;\ell][4\;\ell][a\;i][\sigma_{\ell,4}] \\
{} &=& [\sigma_{\ell,4}][4\;a][\ell\;i][\sigma_{a,4}][a\;\ell] \\
{} &=& [\sigma_{\ell,4}][4\;a][\ell\;i][\sigma_{i,4}][\sigma_{a,4}][\sigma_{i,4}][a\;\ell] \\
{} &=& [\sigma_{\ell,4}][4\;a][\ell\;i][\sigma_{\ell,4}][\sigma_{i,4}][a\;\ell]= \widetilde{\alpha}([\sigma_{i,\ell}])\widetilde{\alpha}([i\;4]).
\end{array}$$

Si $\{i,4\} \cap \{k,\ell\}=\{\ell\}=\{i\}$, alors $\{a,b\}=\{a,4\}$, avec $a \neq k$ et $a \neq i$ et $\tau(k)=k$ et $\tau(\ell)=4$. Donc

$$\begin{array}{ccl}
\widetilde{\alpha}([i\;4])\widetilde{\alpha}([\sigma_{k,i}])&=&[a\;k][\sigma_{i,4}][\sigma_{i,4}][k\;i][a\;4][\sigma_{i,4}] \\
{} &=& [a\;k][k\;i][a\;4][\sigma_{i,4}] \\
{} &=& [a\;i][k\;4][a\;k][\sigma_{i,4}]=\widetilde{\alpha}([\sigma_{k,4}])\widetilde{\alpha}([i\;4]).
\end{array}$$

Si $\{i,4\} \cap \{k,\ell\}=\{i,4\}$, alors $k=4$ et $\{a,b\} \cap \{i,4\}=\varnothing$ et $\tau(k)=i$ et $\tau(\ell)=4$. Donc

$$\begin{array}{ccl}
\widetilde{\alpha}([i\;4])\widetilde{\alpha}([\sigma_{4,i}])&=&[a\;b][\sigma_{i,4}][\sigma_{i,4}][4\;i][a\;b][\sigma_{i,4}] \\
{} &=& [a\;b][4\;i][a\;b][\sigma_{i,4}] \\
{} &=& [4\;i][a\;b][a\;b][\sigma_{i,4}]=\widetilde{\alpha}([\sigma_{i,4}])\widetilde{\alpha}([i\;4]).
\end{array}$$

Donc $\widetilde{\alpha}$ préserve toutes les relations données dans la proposition~\ref{présentation AutWn}. De ce fait, $\widetilde{\alpha}$ se prolonge en un morphisme de groupes de $\Out(W_4)$ dans lui-même. Ceci conclut la démonstration.
\hfill\qedsymbol

\bibliographystyle{alphanumfr}
\bibliography{bibliographie}
\noindent \begin{tabular}{l}
Laboratoire de mathématique d'Orsay\\
UMR 8628 CNRS Univ. Paris-Sud\\
Université Paris-Saclay\\
92405 ORSAY Cedex, FRANCE\\
{\it e-mail: yassine.guerch@math.u-psud.fr}
\end{tabular}
\end{document}